\documentclass[11pt]{article}
\usepackage{amsfonts}
\usepackage{amsmath}
\usepackage{epsfig}

\usepackage{amsfonts}
\usepackage{epsfig}

\newcommand{\dd}{\;\mathrm{d}}

\newenvironment{ack}{\noindent{\small{A}\tiny{CKNOWLEDGMENTS}.}}{} 
\newenvironment{keyword}{\noindent{\small{K}\tiny{EYWORDS}.}}{} 
\newenvironment{pf*}[1]{\noindent {\it #1}}{$\Box$} 
\newenvironment{pf}{\noindent {\it Proof.}}{$\Box$} 
\newcommand{\e}{\mathrm{e}}
\newcommand{\Rset}{\mathbb{R}} 
\newcommand{\Cset}{\mathbb{C}} 
\newcommand{\Qset}{\mathbb{Q}} 
\newcommand{\Zset}{\mathbb{Z}} 
\newtheorem{thm}{Theorem} 
\newtheorem{defn}[thm]{Definition} 
\newtheorem{prop}[thm]{Proposition}

\newcommand{\PP}{\mathbb{P}} 
\newcommand{\TT}{\mathbb{T}}
\newcommand{\Li}{\mathrm{Li}} 
\newcommand{\re}{\mathop{\mathrm{Re}}} 
\newcommand{\im}{\mathop{\mathrm{Im}}} 
\newcommand{\ii}{\mathrm{i}} 
\newcommand{\Lf}{\mathrm{L}} 
\newenvironment{classification}{\noindent{\bf Classification $\,$   }}{}

\begin{document}

\title{Mahler measures and computations with regulators}




\author{Matilde N. Lal\'{\i}n\footnote{This material is partially based upon work supported by the National Science Foundation under agreement No. DMS-0111298.}\\ {\small \it Department of Mathematics, University of British Columbia } \\{ \small \it  Vancouver, BC V6T-1Z2, Canada}\\ {\small \texttt{mlalin@math.ubc.ca}}}

\maketitle 

\begin{abstract}
In this work we apply the techniques that were developed in \cite{L3} in order to study several examples of multivariable polynomials whose Mahler measure is expressed in terms of special values of the Riemann zeta function or Dirichlet L-series. The examples may be understood in terms of evaluations of regulators. Moreover, we apply the same techniques to the computation of generalized Mahler measures, in the sense of Gon and Oyanagi \cite{GO}.
\end{abstract}



\section{Introduction}

The (logarithmic) Mahler measure of a Laurent polynomial $P \in \Cset[x_1^{\pm1}, \dots, x_n^{\pm1}]$ is defined by
\begin{eqnarray}\label{defi}
m(P)& := & \int_0^1 \dots \int_0^1 \log |P(\e^{2 \pi \ii \theta_1}, \dots , \e^{2 \pi \ii \theta_n}) | \dd \theta_1 \dots \dd \theta_n,\\
&=&\frac{1}{(2\pi)^n} \int_{\TT^n} \log|P(x_1, \dots, x_n)|\frac {\dd x_1}{x_1} \dots \frac{\dd x_n}{x_n}.
\end{eqnarray}
where $\TT^n=\{|x_1|=\dots = |x_n|=1\}$. 

The Mahler measure of a one-variable polynomial has a simple expression in terms of the roots, due to Jensen's formula. The several-variable case, however, is much more complicated to describe.

Many formulas for specific examples  have been computed, especially for polynomials in two and three variables. An interesting fact is that many of these formulas express the Mahler measure of a polynomial in terms of special values of the Riemann zeta function or Dirichlet L-series. A typical example is Smyth's formula \cite{S1}
\begin{equation} \label{eq:smy}
m(1+x+y+z)= \frac{7}{2\pi^2} \zeta(3).
\end{equation}
Moreover, L-series of varieties also appear in these kinds of formulas.

Inspired by these results, Deninger \cite{D} established the relation between Mahler measure and regulators in many cases. More precisely, he wrote for $P \in \Cset[x_1, \dots, x_n]$,
\begin{equation}\label{eq:deninger}
m(P)  = m(P^*) + \frac{1}{(-2\pi\ii)^{n-1}} \int_\Gamma \eta_n(n)(x_1,\dots,x_n),
\end{equation}
where
\[\Gamma= \{P(x_1,\dots, x_n) =0\} \cap \{|x_1|=\dots=|x_{n-1}|=1, |x_n|\geq 1 \},\]
and $\eta_n(n)$ is Goncharov's regulator, a $C^\infty$ form on the generic point of $\Gamma$ (see Definition \ref{defn:eta}).

Boyd \cite{B2} computed many more numerical examples that fit into this context. Rodriguez-Villegas \cite{RV} applied the ideas of Deninger and developed others to understand and prove many cases in two-variables by making explicit computations with the regulator.  This work was, in certain sense, continued by Boyd and Rodriguez-Villegas \cite{BRV1,BRV2} for two-variable polynomials. One of the ideas in these works is to identify the cases where $\eta_2(2)(x_1,x_2)$ is exact and proceed to the computation of the Mahler measure by means of Stokes's Theorem, obtaining special values of the Bloch-Wigner dilogarithm.

Later Maillot suggested a way to continue these ideas for more variables keeping in mind the cohomological interpretation.  

In \cite{L3} we developed  these ideas  for some three-variable cases involving evaluations of trilogarithms and shed some light on how these computations could be carried on for more variables.
Let us consider the three-variable case. In favorable cases  there is a primitive for $\eta_3(3)(x_1,x_2,x_3)$ (say $\omega$), one may apply Stokes's Theorem and obtain an  integral
 $\int_{\partial \Gamma}\omega$. Maillot proposed a 
way that desingularizes $\partial \Gamma$  to certain $\widetilde{\partial \Gamma}$ that has an algebraic description and nontrivial boundary. This process also allows the possibility that $\omega$ (restricted to $\widetilde{\partial \Gamma}$ ) is exact. If that is the case, we may continue the integration by applying Stokes's Theorem again.

In this work we explore the techniques of \cite{L3} in application to several concrete examples. Most of the examples that we study here are in three variables, but we also include an example in four variables, and we apply the method to compute some generalized Mahler measures.

More precisely, we recover formula (\ref{eq:smy}) as well as  the following results:

\bigskip
\renewcommand{\arraystretch}{2.2}

\begin{tabular}{ll}
$m(1+x+y^{-1}+(1+x+y)z) = \frac{14}{3\pi^2} \zeta(3)$ & Smyth \cite{S2}\\
$m((1-x)(1-y)+(1+x)(1+y)z) = \frac{7}{\pi^2} \zeta(3)$ & Lal\'{\i}n \cite{L}\\
$m(1+x+2y+(1-x)z) = \frac{7}{2\pi^2} \zeta(3) + \frac{ \log 2}{2}$ & Lal\'{\i}n \cite{L}\\
$m((1-y)(1+x)+(1-x)z) = \frac{28}{5\pi^2} \zeta(3)$& Condon \cite{C}\\
$m((1+x_1)(1+x)+(1-x_1)(1+y)z) = \frac{24}{\pi^3} \Lf(\chi_{-4},4) $ & Lal\'{\i}n  \cite{L,L4}\\
\end{tabular}

\bigskip
\renewcommand{\arraystretch}{1}

Generalized Mahler measures were introduced by Gon and Oyanagi \cite{GO} who studied their basic properties, computed some examples, and related them to multiple sine functions and special values of Dirichlet L-functions.

Given $f_1, \dots, f_r \in \Cset[x_1^{\pm 1}, \dots, x_n^{\pm 1}]$, the generalized (logarithmic) Mahler measure of $f_1, \dots, f_r$ is defined by
\begin{equation}\label{defn2}
m(f_1,\dots, f_r) = \frac{1}{(2\pi \ii)^n} \int_{\TT^n} \max\{\log|f_1|, \dots, \log|f_r|\} \frac{\dd x_1}{x_1} \dots  \frac{\dd x_n}{x_n}.
\end{equation}

Notice that for $r=1$ one obtains the classical Mahler measure of $f_1$. On the other hand, for $r=2$, 
\[m(f_1,f_2) = m(f_1+zf_2),\]
where $z$ is a variable that is independent of $x_1, \dots, x_n$.

In this work we use regulators to study  generalized Mahler measures when $f_i=P(x_i)$ for a fixed polynomial (or rational function) $P$. We obtain explicit formulas for
\[m\left( (1-x_1), \dots, (1-x_{n})\right), \quad m\left( \frac{1-x_1}{1+x_1}, \dots, \frac{1-x_{n}}{1+x_{n}}\right), \]
and
\[m\left(1+x_1-x_1^{-1}, \dots, 1+x_{n}-x_{n}^{-1}\right).\]
The first case was known to Gon and Oyanagi \cite{GO}. To our knowledge, the other two cases are new results.

Finally, we study the behavior of the generalized Mahler measure in the general case of $f_i=P(x_{i_1}, \dots, x_{i_n})$ (again, $P$ is fixed) when the number of functions $f_i$ goes to infinity. The conclusion is that the generalized Mahler measure approaches the sup norm of $P$ in $\TT^n$. This is the content of Proposition \ref{prop}. This result may be combined with the aforementioned examples in order to compute limits of sums involving zeta values, such as
\[\lim_{m \rightarrow \infty} \sum_{j=1}^{m-1} (-1)^{j}\binom{2m-1}{2j} \frac{(2j)!(1-2^{2j})}{(2\pi)^{2j}} \zeta(2j+1) = \log 2.\]

\section{A construction for regulators}
Our goal is to apply the techniques of \cite{L3}. In this section we describe our main ingredients. We will be following Goncharov's construction of the regulator on polylogarithmic motivic complexes \cite{G4,G5}.

First recall Zagier's modification of the polylogarithm \cite{Z2}:
\begin{equation}
 \widehat{\mathcal{L}}_n(x) :=  \widehat{\re}_n \left(\sum_{j=0}^{n-1} \frac{2^j B_j}{j!} (\log|x|)^j \Li_{n-j}(x) \right),
\end{equation}
where $B_j$ is the $j$th Bernoulli number, $\Li_k$ is the classical polylogarithm and $\widehat{\re}_k$ denotes $\re$ or $\ii \im$ depending on whether $n$ is odd or even.   For the record, $\mathcal{L}_n$ is defined using $\re_n$ instead of $\widehat{\re}_n$, where $\re_n$ denotes $\re$ or $\im$ depending on whether $n$ is odd or even.

Polylogarithms satisfy functional equations such as
\[\mathcal{L}_n\left( \frac{1}{x} \right) = (-1)^{n-1} \mathcal{L}_n(x)\]
for $n>1$, and 
 \[ \mathcal{L}_n(\bar{x}) = (-1)^{n-1} \mathcal{L}_n(x),\quad \mathcal{L}_n(x) + \mathcal{L}_n(-x) = \frac{\mathcal{L}_n(x^2)}{2^{n-1}}.\]

For $n=2$,  one obtains the Bloch Wigner dilogarithm,
\begin{equation}
\mathcal{L}_2(x)=D(x) =  \im(\Li_2(x) ) + \log |x| \arg(1-x),
\end{equation}
which satisfies the well-known five-term relation
\begin{equation}
D(x) + D(1-xy) + D(y) + D\left( \frac{1-y}{1-xy} \right) + D\left( \frac{1-x}{1-xy} \right) = 0.
\end{equation}
In particular, $D(x)=-D(1-x)$. A useful functional equation for $n=3$ is
\begin{equation}
\mathcal{L}_3(x) + \mathcal{L}_3(1-x) +\mathcal{L}_3\left(1-\frac{1}{x}\right) = \zeta(3).
\end{equation}

Finally, for an element $z$ of a field $F$, we denote by $\{z\}_n$ the class of $z$ in $\mathcal{B}_n(F) := \Zset[\PP^1_F]/\mathcal{R}_n(F)$, where $\mathcal{R}_n(F)$ is a certain subgroup of $\Zset[\PP^1_F]$ (see definition in \cite{G3}, for instance) that conjecturally consists of all the rational functional equations of the $n$-polylogarithm in $F$.

Given a complex variety $X$, Goncharov considers certain $C^\infty$ forms on the generic point of  $X$,
\begin{equation}
\eta_n(n): \, \bigwedge^n(\Cset(X)^* )_\Qset  \rightarrow \Omega_{X^\infty}^{n-1}(\eta_X),
\end{equation}
\begin{equation}
\eta_{n}(l): \, \mathcal{B}_{n-l+1}(\Cset(X))\otimes_\Qset \bigwedge^{l-1}(\Cset(X)^* )_\Qset  \rightarrow \Omega_{X^\infty}^{l-1}(\eta_X), \qquad l<n.
\end{equation}

The construction is as follows
\begin{defn} \label{defn:eta} Let $x_i$ be rational functions on  $X$.
\[\eta_{n}(n) : x_1 \wedge \dots \wedge x_n \rightarrow\]
\begin{equation}
\mathrm{Alt}_n \left( \sum_{p\geq0}\frac{\log|x_1|}{(2p+1)!(n-2p-1)!} \bigwedge_{j=2}^{2p+1} \dd \log |x_j| \wedge \bigwedge_{j=2p+2}^n  \dd \ii \arg x_j \right),
\end{equation}
where\[ \mathrm{Alt}_m F(t_1,\dots t_m) : =\sum_{\sigma \in S_m} (-1)^{|\sigma|} F(t_{\sigma(1)}, \dots,  t_{\sigma(m)}).\]
\end{defn}

This form has singularities in the zeros and poles of $x_i$. This will not be a problem in the present work, but we refer the reader to \cite{G4} for a detailed discussion. Notice that
\[ \dd \eta_n(n)(x_1,\dots,x_n)  = \widehat{{\re}}_n \left( \frac{\dd x_1}{x_1} \wedge \dots \wedge \frac{\dd x_n}{x_n} \right).\]

We need more notation for the construction of $\eta_{n}(l)$.

For any integers $p\geq 1$ and $k\geq 0$, define
\[\beta_{k,p} : = (-1)^p\frac{(p-1)!}{(k+p+1)!} \sum_{j=0}^{\left[\frac{p-1}{2}\right]} \binom{k+p+1}{2j+1} 2^{k+p-2j}B_{k+p-2j}.\]

Then let
\[\widehat{\mathcal{L}}_{p,q}(x) : = \widehat{\mathcal{L}}_p(x) \log^{q-1}|x| \dd \log|x|, \qquad p\geq 2, \]
\[ \widehat{\mathcal{L}}_{1,q}(x) := (\log|x|\dd \log|1-x|- \log|1-x| \dd \log|x|) \log^{q-1}|x|.\]

\begin{defn}
We have
\[\eta_{n}(l) : \{x\}_{n-l+1} \otimes x_1 \wedge \dots \wedge x_{l-1} \rightarrow\]
\[ \widehat{\mathcal{L}}_{n-l+1}(x) \mathrm{Alt}_{l-1} \left(\sum_{p\geq0}\frac{1}{(2p+1)!(l-1-2p)!} \bigwedge_{j=1}^{2p} \dd \log |x_j| \wedge \bigwedge_{j=2p+1}^{l-1}  \dd \ii \arg x_j \right) \]
\begin{equation}
+ \sum_{1 \leq k,\, 1\leq p \leq l-1} \beta_{k,p} \widehat{\mathcal{L}}_{n-l+1-k,k}(x) \wedge \mathrm{Alt}_{l-1} \left ( \frac{\log|x_1|}{(p-1)!(l-1-p)!}\bigwedge_{j=2}^{p} \dd \log |x_j| \wedge \bigwedge_{j=p+1}^{l-1} \dd \ii\arg x_j   \right).
\end{equation}

\end{defn}
Notice that
\[\eta_n(1) (x)=\widehat{\mathcal{L}}_n(x),\]
\begin{eqnarray*}
\eta_{n}(n)(x_1, 1-x_1, x_2, \dots, x_{n-1}) & = &\dd \eta_{n}(n-1)(x_1, x_2 \dots, x_{n-1} ),\\\\
\eta_{n}(l)(x_1, x_1, x_2, \dots, x_{l-1}) & = & \dd \eta_{n}(l-1)(x_1, x_2 \dots, x_{l-1} ), \qquad l <n.
\end{eqnarray*}
This equality is to be understood in any open set $U \subset X$ where the forms are defined, i.e., excluding the zeros and poles of $x_i$.

To be concrete, we describe the two-variable case. We have
\begin{equation}
\eta(x,y):= \eta_2(2) (x,y) = \log|x| \dd \ii \arg(y) - \log|y| \dd \ii \arg(x).  
\end{equation}
Moreover,
\[\eta_2(2) (x,1-x) = \dd \ii D(x). \]

The forms for $n=3$ are
\[\eta(x,y,z):= \eta_3(3) (x,y,z) =\log|x|\left( \frac{1}{3} \dd \log |y|\,\wedge \dd \log |z| +\dd \ii \arg y\,\wedge \dd \ii \arg z\right )\]
\[ + \log|y|\left( \frac{1}{3} \dd \log |z|\,\wedge \dd \log |x| + \dd \ii \arg z \,\wedge \dd \ii \arg x\right)\]  \[+ \log|z|\left( \frac{1}{3} \dd \log |x|\,\wedge \dd \log |y| + \dd \ii \arg x\, \wedge\dd \ii \arg y\right),\]

\[ \eta_3(3)(x,1-x,y) = \dd \eta_3(2) (x,y),\]

\[\omega(x,y):= \eta_3(2)(x,y)= \ii D(x) \dd \ii \arg y - \frac{1}{3}(\log |x| \dd \log |1-x|- \log |1-x| 
\dd \log |x| )\log |y|  .\]
For future reference, observe that $\eta(x,y,z)$ changes sign under complex conjugation of $x,y,$ and $z$, and that $\omega(x,y)$ is invariant under complex conjugation of both $x$ and $y$.

\section{The relation with Mahler measure}

In this section we summarize the method described in \cite{L3}.
Let $P \in \Cset[x_1^\pm, \dots, x_n^\pm ]$. We may assume, without loss of generality, that $P$ is actually a polynomial. Then we may write
\[P(x_1, \dots, x_n) = a_d(x_1, \dots, x_{n-1})x_n^d + \dots + a_0(x_1, \dots, x_{n-1}).\]
Thus, formally,
\[P(x_1, \dots, x_{n}) = a_d(x_1, \dots, x_{n-1}) \prod_{j=1}^d (x_n-\alpha_j(x_1, \dots, x_{n-1})).\]
By Jensen's formula,
\[m(P) = m(a_d) + \frac{1}{(2 \pi \ii)^{n-1}} \sum_{j=1}^d \int_{\TT^{n-1}} \log^+ | \alpha_j(x_1, \dots, x_{n-1})|
\frac{\dd x_1}{x_1}\wedge \dots \wedge \frac{\dd x_{n-1}}{x_{n-1}}.\]
From Definition \ref{defi},
\begin{equation}\label{eq:Jensen}
m(P)= m(a_d) + \frac{1}{(-2 \pi \ii)^{n-1} } \int_\Gamma \eta_n(n)(x_1, \dots, x_n), 
\end{equation}
where
\[\Gamma= \{P(x_1,\dots, x_n) =0\} \cap \{|x_1|=\dots=|x_{n-1}|=1, |x_n|\geq 1 \}.\]
Thus we have recovered Deninger's expression (\ref{eq:deninger}).

We have to compute the integral. Our strategy is as follows. For a given polynomial, we try to prove that $\eta_n(n)(x_1, \dots, x_n)$ is exact. In order to achieve that, we try to express $x_1 \wedge \dots \wedge x_n$ as a linear combination of elements of the form $y_1 \wedge (1-y_1) \wedge y_2\wedge \dots \wedge y_{n-1}$. At the same time we compute the boundary of $\Gamma$. Then we apply Stokes's Theorem in order to obtain a linear combination of integrals of the form 
\[ \int_{\partial \Gamma} \eta_{n}(n-1)(y_1,\dots, y_{n-1}).\] 
We continue by examining the restriction of  $\{y_1\}_2\otimes y_2 \wedge 
\dots \wedge y_{n-1}$ to each component of $|\partial \Gamma|$ and trying 
to express it as a linear combination of elements of the form 
$\{z_1\}_2\otimes z_1 \wedge z_2 \wedge \dots \wedge z_{n-2}$. This will 
be different for different components, since $\eta_n(n-1)$ is not even closed on $\partial\Gamma$ itself. If we are successful,  then we may obtain a sum of integrals of the form 
\[ \int_{\partial^2 \Gamma} \eta_{n}(n-2)(z_1,\dots, z_{n-2}).\] 
Notice that $\partial^2 \Gamma$ should be empty. However, typically $\partial \Gamma$ will have singularities. Upon desingularizing $\partial \Gamma$, we obtain a set with nontrivial boundary.

If we are fortunate enough, we may continue this process until we compute the integral completely. The success of this method depends on the polynomial $P$. See \cite{L3} for a discussion of sufficient conditions and other details.

\section{Examples}
Now we will apply the machinery that was described in the previous sections in order to understand many examples of Mahler measure formulas in three variables. We will study many of the examples that are known to be related to $\zeta(3)$ and the trilogarithm. 

It is hard in practice to apply the technique when the polynomial has degree higher than 1 in the variable $x_n(=x_3)$.  From now on we will replace $(x_1,x_2,x_3)$  with $(x,y,z)$. Hence, all the examples will be written as
\begin{equation} \label{eq:R}
z=R(x,y)
\end{equation}
where $R(x,y)$ is a rational function with real coefficients. 

In computing the boundary $\partial \Gamma$, we apply the idea proposed by Maillot \cite{L3}. Since
\[ \Gamma= \{z=R(x,y)\} \cap \{|x|=|y|=1, |z| \geq 1\},\]
we have
\[ \partial \Gamma= \{z=R(x,y)\} \cap \{|x|=|y|=|z| =1\}.\]
Notice that $R$ has real coefficients, therefore we can describe
\[\partial \Gamma = \{z-R(x,y)=z^{-1}-R(x^{-1},y^{-1})=0\} \cap \{|x|=|y|=|z|=1\}.\]
Then we may write 
\[\gamma := \partial \Gamma = C \cap \{|x|=|y|=|z|=1\},\]
where $C$ is the curve defined by 
\begin{equation}\label{eq:C}
R(x,y)R(x^{-1},y^{-1})=1.
\end{equation}

For the case in four variables we can follow a similar process.

In most of the cases the denominator of the rational function $R$ is a product of cyclotomic polynomials. We may multiply Eq. (\ref{eq:R}) by the denominator in order to obtain a polynomial equation whose Mahler measure is the same as the Mahler measure of $z-R$. We will call this polynomial $P$. The only exception to this case is $P= 1+x+y^{-1}-(1+x+y)z$, since the measure of $1+x+y$ is not zero. In what follows, we will compute the (logarithmic) Mahler measure of $P$.

\subsection{Smyth's example}

We are going to start with the simplest example in three variables, which was also due to Smyth \cite{B1,S2}:
\begin{equation}
\pi^2 m(1+x+y+z) = \frac{7}{2} \zeta(3) .
\end{equation}

It is easy to see that this problem amounts to computing the Mahler measure of the rational function
\[ z - \frac{1-x}{1-y}. \]

We would like to see that $\eta(x,y,z)$ is exact. The equation with the wedge product yields
\begin{eqnarray*}
 x \wedge y \wedge z & =& x \wedge y \wedge \frac{1-x}{1-y} = x \wedge y \wedge (1-x) - x \wedge y \wedge (1-y)\\
 & = & - x \wedge (1-x) \wedge y -y \wedge (1-y)\wedge x .
 \end{eqnarray*}

In other words,
\[\eta(x,y,z)  = -\eta(x , 1-x, y) - \eta(y,1-y, x).\]

After performing the first integration, we are left to analyze
\[\Delta = -\{x\}_2 \otimes y - \{y\}_2 \otimes x.\]

Eq. (\ref{eq:C}) yields for $|\partial \Gamma|$
\[ (xy-1)(x-y)=0\]
in this case.

When $xy=1$ we obtain
\[\Delta = 2 \{x\}_2 \otimes x.\]

When $x=y$ we obtain
\[\Delta = - 2 \{x\}_2 \otimes x.\]

One could have problems in the cases when $z$ has a pole or is equal to zero. But those correspond to $x=1$ or $y=1$, and $\Delta=0$ in these circumstances.

We obtain
\[ -\omega(x,y) - \omega(y,x) = \pm 2\omega(x,x), \]
which yields
\[m(P) = \frac{1}{4 \pi^2} \int_\gamma 2\omega(x,x),\]
where $\gamma = \partial \Gamma$.

We now need to check the path of integration $\gamma$. Since the equations $xy=1$ and $x=y$ intersect in $(\pm 1 , \pm 1)$, there are four paths, which can be parameterized as
\[\begin{array}{cccc}
 x = \e^{\alpha \ii}, & 0 \leq \alpha \leq \pi,  & y = x, & -2\{x\}_3; \\\\
 x = \e^{\alpha \ii}, & \pi \geq \alpha \geq 0, & yx =1, & 2\{x\}_3;\\\\
 x = \e^{\alpha \ii}, &  0 \geq \alpha \geq -\pi, & y = x, &-2\{x\}_3;\\\\
 x = \e^{\alpha \ii}, &  -\pi \leq \alpha \leq 0, &  yx =1, & 2\{x\}_3,
 \end{array}\]
where the right hand side column indicates the primitives for the pullback of $\Delta$ to each component of $\gamma$.
\begin{figure}
\centering
\includegraphics[width=12pc]{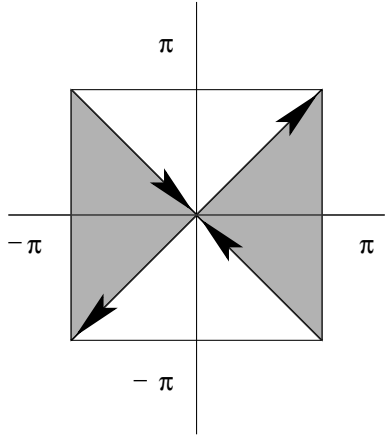}
\caption{\label{mat3smyth} Integration set for $1+x+y+z$. }
\end{figure}
Figure \ref{mat3smyth} shows the integration set for the problem in terms of $\arg x$ and $\arg y$. The shaded region corresponds to the original surface $\Gamma$, where the first integration is performed, and $\gamma$ is the boundary of this surface.

Finally we obtain
\begin{equation}
m(P)= \frac{1}{4 \pi^2}8(\mathcal{L}_3(1) - \mathcal{L}_3(-1)) = \frac{7}{2\pi^2} \zeta(3).
\end{equation}

It is important to notice that the orientation for $\gamma$ is not the one that we could naively expect without taking into account the singularities at $\pm (1,1)$. That is to say, if we think of the diagonals in the picture as two circles, the boundary would be zero and we would get the erroneous result $m(P)=0$. The fact that different paths in each circle are oriented differently is crucial for the result.

The orientation of $\Gamma$ is not as easy to determine. The simplest way to do this is to choose any orientation for $\Gamma$ and use the induced orientation in $\gamma$. If the chosen orientation for $\Gamma$ is incorrect, then we obtain $m(P)$ with the wrong sign, but that can be easily corrected since we know $m(P)$ must be positive.

\subsection{Another example due to Smyth}\label{Smyth2}

We will now be concerned with another example due to Smyth \cite{S2},
\[  P=1+x+y^{-1}-(1+x+y)z.\]

In this case, the equation for the wedge product yields 
 \begin{eqnarray*}
 x \wedge y \wedge z & =& x \wedge y  \wedge (1+x+y^{-1})- x \wedge y \wedge (1+x+y )\\
 & = & - x \wedge y^{-1}  \wedge (1+x+y^{-1})- x \wedge y \wedge (1+x+y ).
\end{eqnarray*}
But
\begin{eqnarray*}
x \wedge y \wedge (1+x+y ) &=& \frac{x}{y} \wedge y \wedge (1+x+y )\\
&= & \frac{x}{y} \wedge (x+y) \wedge (1+x+y) - \frac{x}{y} \wedge \left(1+\frac{x}{y}\right) \wedge (1+x+y) \\
& = & (-x-y) \wedge (1+x+y) \wedge  \frac{x}{y} - \left(- \frac{x}{y}\right) \wedge \left(1+\frac{x}{y}\right) \wedge (1+x+y).
\end{eqnarray*}


Then  we need to analyze
\[ \Delta = -\left\{-x-\frac{1}{y}\right\}_2 \otimes xy + \left\{ -xy \right\}_2 \otimes\left (1+x + \frac{1}{y}\right) - \{-x-y\}_2 \otimes \frac{x}{y} + \left\{ -\frac{x}{y} \right\}_2 \otimes (1+x + y). \]
The technique of Eq. (\ref{eq:C}) yields for $|\partial \Gamma|$
\[ (x-x^{-1})(y-y^{-1}) = 0. \]

If $y=-1$,
\[ \Delta =- \left\{1-x\right\}_2 \otimes(- x) + \left\{ x \right\}_2 \otimes x  - \{1-x\}_2 \otimes (-x) + \left\{ x \right\}_2 \otimes x \]
\[ = 4\{x\}_2 \otimes x . \]

If $y= 1$,
\[ \Delta =- 2 \left\{-1-x\right\}_2 \otimes x +2 \left\{ -x \right\}_2 \otimes \left (2+x \right) =2 \left\{2+x\right\}_2 \otimes x +2 \left\{ -x \right\}_2 \otimes \left (2+x \right).  \]

We will use the five-term relation starting with $\{2+x\}_2$ and $\{-x\}_2$,
\[2\{2+x\}_2+2\{-x\}_2 + \{(1+x)^2\}_2 = 0.\]

We obtain (by using $-\{a\}_2=\{1-a\}_2$),
\[ \Delta = -2\{-x\}_2\otimes x - \{(1+x)^2\}_2 \otimes x - 2\{2+x\}_2 \otimes(2+x) - \{(1+x)^2\}_2 \otimes (2+x)  \]
\[ = -2\{-x\}_2\otimes (-x) - 2\{2+x\}_2 \otimes(2+x) + \{-2x-x^2\}_2 \otimes (-2x-x^2).  \]

If $x=-1$,
\[ \Delta =- \left\{1-\frac{1}{y}\right\}_2 \otimes y - \left\{ y \right\}_2 \otimes y + \{1-y\}_2 \otimes (-y) + \left\{ \frac{1}{y} \right\}_2 \otimes y  \]
\[=-4\{y\}_2 \otimes y. \]

If $x=1$,
\[ \Delta =- \left\{-1-\frac{1}{y}\right\}_2 \otimes y +\left\{ -y \right\}_2 \otimes\left (2 + \frac{1}{y}\right) - \{-1-y\}_2 \otimes \frac{1}{y} + \left\{ -\frac{1}{y} \right\}_2 \otimes (2+ y). \]

We will use the five-term relation starting with $\left\{\frac{1}{y}\right\}_2$ and $\{-1-y\}_2$,
\[ 2\left\{ -\frac{1}{y} \right\}_2 + 2\{ -1-y\}_2  + \{-2y-y^2\}_2 = 0 \]
then
\[ 2 \left\{ -\frac{1}{y} \right\}_2 + 2 \{ -1-y\}_2 = \{(1+y)^2\}_2.\]

Now
\[ - \{-1-y\}_2 \otimes \frac{1}{y} + \left\{ -\frac{1}{y} \right\}_2 \otimes (2+ y) \]
\[= -\frac{1}{2} \{(1+y)^2\}_2 \otimes \frac{1}{y} + \left\{ -\frac{1}{y} \right\}_2 \otimes \frac{1}{y} + \frac{1}{2} \{(1+y)^2\}_2 \otimes (2+y) -  \{ -1 -y\}_2 \otimes (2+y) \]
\[= -\frac{1}{2} \{-2y-y^2\}_2 \otimes (-2y-y^2) + \left\{ -\frac{1}{y} \right\}_2 \otimes \left(-\frac{1}{y}\right) + \{ 2+y\}_2 \otimes (2+y).\]

Then we obtain
\[\Delta =- \frac{1}{2} \left\{-\frac{2}{y}-\frac{1}{y^2}\right\}_2 \otimes \left(-\frac{2}{y}-\frac{1}{y^2}\right) + \left\{ -y \right\}_2 \otimes \left(-y\right) + \left\{ 2+\frac{1}{y}\right\}_2 \otimes \left(2+\frac{1}{y}\right) \]
\[-\frac{1}{2} \{-2y-y^2\}_2 \otimes (-2y-y^2) + \left\{ -\frac{1}{y} \right\}_2 \otimes \left(-\frac{1}{y}\right) + \{ 2+y\}_2 \otimes (2+y) .\]

We may need to take into account the poles or zeros of $z$. But those are at the points $(x,y) = (\zeta_6,\zeta_6^{-1})$, and they do not affect the integration because they are a set of points and the integration is in dimension 2.

\begin{figure}
\centering
\includegraphics[width=12pc]{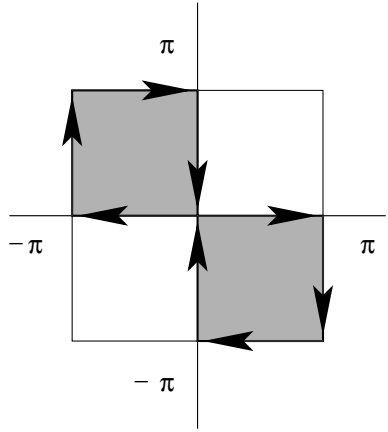}
\caption{\label{mat3smyth2} Integration set for $1+x+y^{-1}-(1+x+y)z$.}
\end{figure}

We compute the boundary $\gamma$ (see Fig. \ref{mat3smyth2}) and primitives for pullbacks of $\Delta$.

\[\begin{array}{cccc}
 x = \e^{\alpha \ii}, & -\pi \leq \alpha \leq 0,  & y = -1, & 4\{x\}_3; \\\\
 x = \e^{\alpha \ii}, &  0 \geq \alpha \geq -\pi, & y = 1, & - 2\{-x\}_3 -2 \{2+x\}_3+\{-2x-x^2\}_3;\\\\
 y = \e^{\alpha \ii}, &  0 \leq \alpha \leq \pi,  & x = -1, &-4\{y\}_3; \\\\
 y = \e^{\alpha \ii}, &  \pi \geq \alpha \geq 0,  &  x = 1, & -\frac{1}{2}\{-2y^{-1}-y^{-2}\}_3 + \{-y\}_3\\
 &&&+\{2+y^{-1}\}_3 -\frac{1}{2}\{-2y-y^{2}\}_3 +\{-y^{-1}\}_3+\{2+y\}_3;
\end{array}\]
\[\begin{array}{cccc}
 x = \e^{\alpha \ii}, & \pi \geq \alpha \geq 0,  & y = -1, & 4\{x\}_3; \\\\
 x = \e^{\alpha \ii}, &  0 \leq \alpha \leq \pi, & y = 1, & - 2\{-x\}_3 -2 \{2+x\}_3+\{-2x-x^2\}_3;\\\\
 y = \e^{\alpha \ii}, &  0 \geq \alpha \geq -\pi,  & x = -1, &-4\{y\}_3; \\\\
 y = \e^{\alpha \ii}, &  -\pi \leq \alpha \leq 0,  &  x = 1, & -\frac{1}{2}\{-2y^{-1}-y^{-2}\}_3 +\{-y\}_3\\
 &&&+\{2+y^{-1}\}_3 -\frac{1}{2}\{-2y-y^{2}\}_3+ \{-y^{-1}\}_3+\{2+y\}_3.
 \end{array}\]

Then we obtain
\[4\pi^2m(P) = 16 (\mathcal{L}_3(1) - \mathcal{L}_3(-1) ) + 4( 2\mathcal{L}_3(-1)-3\mathcal{L}_3(1)  +2 \mathcal{L}_3(3) - \mathcal{L}_3(-3)) \]
\[ = 4 \mathcal{L}_3(1) - 8 \mathcal{L}_3(-1) + 8 \mathcal{L}_3(3) - 4 \mathcal{L}_3(-3).  \]

It will be necessary to use the identity:
\begin{equation}\label{eq:22}
2\mathcal{L}_3(3) - \mathcal{L}_3(-3) =  \frac{13}{6} \zeta(3) 
\end{equation}
which is essentially Lemma 6 in Smyth's \cite{S2}.

Finally,
\begin{equation}
 m(P) = \frac{14}{3\pi^2} \zeta(3).
\end{equation}

It is worth noting that this example involves the evaluation of the element $2\{3\}_3-\{-3\}_3$. As we noted above, Smyth \cite{S2} also encountered this difficulty in his direct computation. 

In order to evaluate the element $2\{3\}_3-\{-3\}_3$ in $\mathcal{B}_3(\Qset)$, we need to use Goncharov's 22-term relation (see \cite{G3}). We follow the notation in Zhao \cite{Zh}, and take the relation with $(a,b,c)=(3,-1,1)$, obtaining
\[4 \{3\}_3+2 \left\{\frac{1}{3}\right\}_3-3\left\{-\frac{1}{3}\right\}_3+6\{-1\}_3-2 \{1\}_3=0.\]
Using that $\{-1\}_3=-\frac{3}{4}\{1\}_3$ and $\left\{\frac{1}{x}\right\}_3=\{x\}_3$,
\[6 \{3\}_3-3\left\{-3\right\}_3-\frac{13}{2}\{1\}_3=0,\]
thus, proving Eq. (\ref{eq:22}). We owe this argument to the referee of this paper.

\subsection{Another three-variable example \label{sec15}}
The following example was first computed in \cite{L}. It is easier to consider the following rational function
\[ z - \frac{(1-x)(1-y)}{(1+x)(1+y)}. \]

For the wedge product we have,
\begin{eqnarray*}
 x\wedge y \wedge z & = & x\wedge y \wedge (1-x) +  x\wedge y \wedge (1-y) -  x\wedge y \wedge (1+x) -  x\wedge y \wedge (1+y)\\
 & = &- x\wedge (1-x) \wedge y+   y \wedge (1-y)\wedge x +  (-x)\wedge (1+x)\wedge y -  (-y) \wedge (1+y) \wedge x.
 \end{eqnarray*}

Thus, we need to consider
\[ \Delta = -\{x\}_2 \otimes y +\{y\}_2 \otimes x + \{-x\}_2 \otimes y -   \{-y\}_2 \otimes x. \]
This time Eq. (\ref{eq:C}) implies
\[(xy+1) (x+y)= 0.\]

When $xy=-1$, 
\[\Delta= 2 \{x\}_2\otimes x - 2 \{-x\}_2 \otimes (-x). \]

When $x=-y$,
\[\Delta= -2 \{x\}_2\otimes x + 2 \{-x\}_2 \otimes (-x). \]

The poles or zeros in this case occur with $x=\pm1$ and $y=\pm1$, but we always obtain $\Delta =0$ and they do not affect the integration.

We now need to check the integration path $\gamma$ and primitives for pullbacks of $\Delta$.

\begin{figure} 
\centering
\includegraphics[width=12pc]{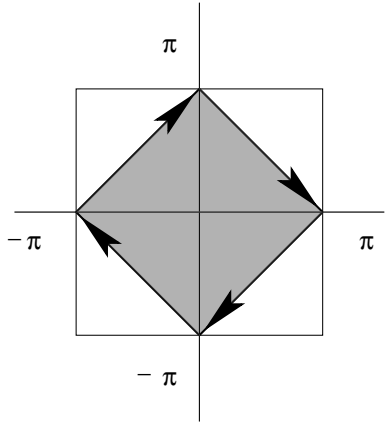}
\caption{\label{mat3l} Integration set for $(1-x)(1-y)-(1+x)(1+y)z$.}
\end{figure}

\[\begin{array}{cccc} 

 x = \e^{\alpha \ii}, & 0 \geq \alpha \geq -\pi,  & y = -x^{-1}, & 2\{x\}_3 -2 \{-x\}_3; \\\\
 x = \e^{\alpha \ii}, &  \pi \geq \alpha \geq 0, & y = -x, & -2\{x\}_3 +2 \{-x\}_3;\\\\
 x = \e^{\alpha \ii}, &  0 \leq \alpha \leq \pi,  &y = -x^{-1}, & 2\{x\}_3 -2 \{-x\}_3;\\\\
 x = \e^{\alpha \ii}, &  -\pi \leq \alpha \leq 0, &y = -x,   & -2\{x\}_3 +2 \{-x\}_3.
 \end{array}\]

Therefore, we obtain,
\[4 \pi^2 m(P) = 16(\mathcal{L}_3(1)-\mathcal{L}_3(-1)), \]
\begin{equation} \label{eq:mine}
 m(P)=  \frac{7}{\pi^2} \zeta(3).
\end{equation}

\subsection{An example with non-trivial symbol}
Now we will study an example that is of different nature because its symbol in $K$-theory is not trivial (see \cite{L3}). In other words, $\eta(x,y,z)$ is not exact in this case. This example was first computed in \cite{L}. As before, we will consider a simpler form than the one in \cite{L}, that is to say the rational function \[  z - \frac{1+x+2xy}{1-x}. \]

We have
\[ x \wedge y \wedge z = x \wedge y \wedge (1+x+2xy) - x \wedge y \wedge (1-x).\]
Now we use that
\[ x \wedge 2y \wedge (1+x+2xy) = (-x) \wedge (-2y) \wedge (1+x(1+2y))\]
\[ =(-x(1+2y)) \wedge (-2y) \wedge (1+x(1+2y)) - (1+2y) \wedge (-2y) \wedge (1+x(1+2y)).\]

Then
\begin{eqnarray*} 
 x \wedge y \wedge z &=& 2\wedge x \wedge z +  x \wedge (1-x) \wedge (2y) - (-x(1+2y)) \wedge  (1+x(1+2y))\wedge (-2y) \\
 && + (-2y)\wedge (1+2y) \wedge (1+x(1+2y)).
 \end{eqnarray*}

We need to analyze
\[\Delta = \{x\}_2\otimes (2y) - \{-x(1+2y)\}_2 \otimes (-2y) + \{-2y\}_2 \otimes (1+x(1+2y)),  \]
and then we also need to compute the integral of $\eta(2,x,z)$ (the nontrivial part in $K$-theory, i.e., the non-exact part). 

By Eq. (\ref{eq:C}) on $|\partial \Gamma|$,
\[xy=-1, \qquad x = -1, \qquad \mbox{or} \qquad y = -1.\]

For $x =-1$,
\[\Delta = - \{1+2y\}_2 \otimes (-2y) + \{-2y\}_2 \otimes (-2y) = 2\{-2y\}_2 \otimes (-2y).\]

For $y = -1$, $\Delta = 0$.

For $xy=-1$,
\[ \Delta =\left\{-\frac{1}{y}\right\}_2\otimes (2y)- \left\{\frac{1+2y}{y}\right\}_2 \otimes (-2y) + \{-2y\}_2 \otimes \left(-1-\frac{1}{y}\right).\]

But
\[-\left\{\frac{1+2y}{y}\right\}_2 = \left\{ -1 - \frac{1}{y} \right \}_2,\]
and we may use the five-term relation (ignoring torsion) in order to get
\begin{equation}\label{eqL2}
 \left\{ -1 - \frac{1}{y} \right \}_2 + \{-2y\}_2 + \{ -1 -2 y \}_2 + \left\{ -\frac{1}{y}\right\}_2 = 0. 
 \end{equation}

Then
\[ \Delta = - \{-2y\}_2\otimes (-2y) - \{-1-2y\}_2 \otimes (-2y)  -\left\{ -1 - \frac{1}{y} \right \}_2 \otimes\left(-1-\frac{1}{y}\right) \]
\[ -\{-1-2y\}_2 \otimes \left( -1-\frac{1}{y}\right) - \left\{ -\frac{1}{y}\right\}_2 \otimes  \left(-1-\frac{1}{y}\right)\]
\[  = - \{-2y\}_2\otimes (-2y) -\left\{ -1 - \frac{1}{y} \right \}_2 \otimes\left( -1-\frac{1}{y}\right)+ \left\{ 1 +\frac{1}{y}\right\}_2 \otimes  \left( 1+\frac{1}{y}\right)  \]
\[+\{2+2y\}_2 \otimes (2+2y).\]

There are zeros of $z$ when $1+x+2xy=0$, and that only can happen in the unit torus if $(x,y) = (1,-1)$ and that is just a point. There are poles for $x=1$, in this case
\[ \Delta = -\{-1-2y\}_2\otimes(2y) + \{-2y\}_2\otimes (2+2y).\]

By the five-term relation (\ref{eqL2}),
\[ \Delta = \left\{-1-\frac{1}{y}\right\}_2 \otimes (2y) + \{-2y\}_2\otimes (2y) - \{-y\}_2\otimes(2y) \]
\[- \left\{-1-\frac{1}{y}\right\}_2 \otimes (2+2y) + \{2+2y\}_2\otimes (2+2y) + \{-y\}_2\otimes(2+2y).\]

This integrates to 
\[ \Omega =\{-2y\}_3 + \{2+2y\}_3 +\left\{1+\frac{1}{y}\right\}_3-\left\{-1-\frac{1}{y}\right\}_3, \]
hence to zero when $y$ moves in the unit circle.

We now need to check the integration path $\gamma$ and the primitives for pullbacks of $\Delta$. 
\begin{figure} 
\centering
\includegraphics[width=12pc]{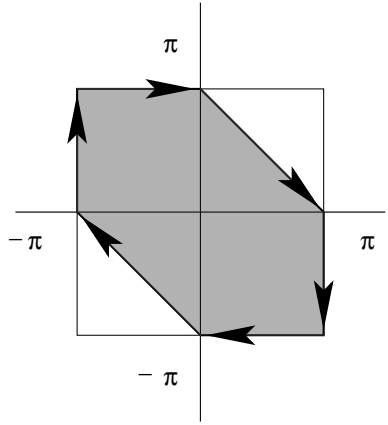}
\caption{\label{mat3l2}Integration set for $1+x+2xy-(1-x)z$.}
\end{figure}

\[\begin{array}{cccc} 
 y = \e^{\alpha \ii}, &   0 \leq \alpha \leq \pi,  &x = -1, & 2\{-2y\}_3; \\\\
 y = \e^{\alpha \ii}, & -\pi \leq \alpha \leq 0,  & x = -y^{-1}, & -\{-2y\}_3 -\{-1-y^{-1}\}_3+\{1+y^{-1}\}_3+\{2+2y\}_3;  \\\\
 y = \e^{\alpha \ii}, &   0 \geq \alpha \geq -\pi,  &x = -1, & 2\{-2y\}_3; \\\\
 y = \e^{\alpha \ii}, & \pi \geq \alpha \geq 0,  & x = -y^{-1}, & -\{-2y\}_3 -\{-1-y^{-1}\}_3+\{1+y^{-1}\}_3+\{2+2y\}_3.  \\
 
 \end{array}\]

We obtain,
\[4\pi^2m(P)= \int \eta(2,x,z)+2(2\mathcal{L}_3(2) - 2\mathcal{L}_3(-2)) +2(\mathcal{L}_3(4) + 2\mathcal{L}_3(2) -2 \mathcal{L}_3(-2)).\]

But $\{4\}_3=4\{2\}_3+4\{-2\}_3$, then
\[4\pi^2m(P)= -\int \eta(2,x,z)+16\mathcal{L}_3(2).\]
Since $\{-1\}_3 + \{2\}_3 +\{2\}_3 = \{1\}_3$, we obtain
\[ \{2\}_3 = \frac{1}{2}(\{1\}_3-\{-1\}_3) = \frac{7}{8} \{1\}_3.\]
Then
\[4\pi^2m(P)= -\int \eta(2,x,z)+14 \zeta(3). \]

We still need to compute
\[ \int_B \eta(2,x,z).\]
where $B=\{(x,y)\in \TT^2\, | \,-\pi\leq \arg x,\arg y, \arg x + \arg y \leq \pi\}$.

Since $z(1-x) = 1+x+2xy$,
\[ \frac{\dd z}{z} = \frac{z +2y +1}{z(1-x)} \dd x + \frac{2x}{z(1-x)} \dd y.\]
But $|x|=1$, so
\[ \dd \arg x \wedge \dd \arg z = -\re \left( \frac{\dd x}{x} \wedge \frac{\dd z}{z} \right) = -\re \left( \frac{\dd x}{x} \wedge \frac{2x}{z(1-x)} \dd y \right).\]

Thus,
\[-\int_B \eta(2,x,z) =  \log 2\int_B \dd \arg x \wedge \dd \arg z =- \re \left( \log 2 \int_B \frac{1}{1+\frac{1+x}{2xy}}\frac{\dd x}{x} \frac{\dd y}{y} \right).\]

We need to consider
\[\int_B \sum_{k=0}^\infty \left(-\frac{1+x}{2xy}\right)^k \frac{\dd x}{x} \frac{\dd y}{y}.\]
Setting $x=\e^{\ii \alpha}$, $y=\e^{\ii \beta}$,
\begin{equation}\label{log2}  
  = -\left( \int_{-\pi}^0 \int_{-\pi-\alpha}^\pi+ \int_{0}^\pi \int_{-\pi}^{\pi-\alpha} \right)\sum_{k=0}^\infty \left(-\frac{(1+\e^{-\ii \alpha})\e^{-\ii\beta}}{2}\right)^k \dd \beta \dd \alpha.
  \end{equation}
  
We consider each term in the series separately,
\[- \int_{-\pi}^0 \int_{-\pi-\alpha}^\pi (-(1+\e^{-\ii \alpha})\e^{-\ii\beta})^k \dd \beta \dd \alpha.  \]
If $k=0$, the above integral is $-\frac{3\pi^2}{2}$. If not,
\[=- \int_{-\pi}^0 (-(1+\e^{-\ii \alpha}))^k \frac{(-1)^k(1- \e^{\ii k \alpha})}{-\ii k} \dd \alpha =  -\frac{\ii}{k}\int_{-\pi}^0 (1+\e^{-\ii \alpha})^k (1- \e^{\ii k \alpha}) \dd \alpha\]
\[=-\frac{\ii}{k} \sum_{l=0}^k \binom{k}{l} \int_{-\pi}^0  \e^{-\ii l \alpha} (1- \e^{\ii k \alpha}) \dd \alpha \]
\[= -\frac{\ii}{k} \left (\frac{2 \ii}{k}(1-(-1)^k) + \sum_{l=1}^{k-1} \binom{k}{l}\ii \left( \frac{1-(-1)^l}{l}+ \frac{1-(-1)^{k-l}}{k-l} \right)\right)\]
\[= \frac{2}{k}  \sum_{l=1}^{k} \binom{k}{l}\frac{1-(-1)^l}{l}.\]

Hence, in order to evaluate the first term in Eq. (\ref{log2}), we need to evaluate
\[  \sum_{k=1}^\infty \frac{2}{k 2^k}  \sum_{l=1}^{k} \binom{k}{l}\frac{1-(-1)^l}{l}= 2 \sum_{l=1}^\infty \frac{1-(-1)^l}{l} \sum_{k=l}^\infty \binom{k}{l}\frac{1}{k 2^k}.\]

But 
\[  \sum_{k=l}^\infty \binom{k}{l}\frac{\lambda^k}{k} = \frac{\lambda^l}{l!}  \sum_{k=l}^\infty (k-1) \dots (k-l+1)\lambda^{k-l} = \frac{\lambda^l}{l!} \frac{\partial^{l-1}}{\partial \lambda^{l-1}}\left( \frac{1}{1-\lambda}\right) \]
\[ =\frac{\lambda^l}{l!} \frac{(l-1)!}{(1-\lambda)^l}= \frac{\lambda^l}{l(1-\lambda)^l}.\]

Using this with $\lambda = \frac{1}{2}$, 
\[  \sum_{k=1}^\infty \frac{2}{k 2^k}  \sum_{l=1}^{k} \binom{k}{l}\frac{1-(-1)^l}{l}=
2 \sum_{l=1}^\infty \frac{1-(-1)^l}{l^2} = 3 \zeta(2) = \frac{\pi^2}{2}.\]

The second integral is the same, so as a conclusion, we get
\[\int_B \sum_{k=0}^\infty \left(-\frac{1+x}{2xy}\right)^k \frac{\dd x}{x} \frac{\dd y}{y} = -2  \pi^2.\]

And, 
\[- \int_B \eta(2,x,z) = 2 \pi^2 \log 2.\]

Finally, the whole formula becomes
\begin{equation}
m(P) = \frac{7}{2\pi^2} \zeta(3) + \frac{\log 2}{2}.
\end{equation}

\subsection{Condon's example}

The last and most complex example that we will analyze in three variables was discovered numerically by Boyd and proved by Condon \cite{C}. It may be expressed in the following way:
\[ z - \frac{(1-y)(1+x)}{1-x}.\]

The wedge product equation becomes:
\begin{eqnarray*}
 x \wedge y \wedge z & = & x \wedge y \wedge (1-y) + x \wedge y \wedge (1+x) - x \wedge y \wedge (1-x) \\
& = & y \wedge (1-y) \wedge x - (-x) \wedge (1+x) \wedge y + x \wedge (1-x) \wedge y.
\end{eqnarray*}

Hence we need to consider
\[ \Delta = \{y\}_2\otimes x - \{-x\}_2 \otimes y + \{x\}_2 \otimes y.\]

Eq. (\ref{eq:C}) implies (on $|\partial \Gamma|$)
 \[ \left( \frac{1+x}{1-x}\right)^2 = \frac{y}{(1-y)^2}. \]

We now use the five term relation,
\[ \{x\}_2 + \{-1\}_2 + \{1+x\}_2 + \left\{ \frac{1-x}{1+x}\right\}_2 + \left\{ \frac{2}{1+x}\right\}_2 = 0, \]
\[ \{x\}_2 - \{-x\}_2 + \left\{ \frac{1-x}{1+x}\right\}_2 - \left\{ \frac{x-1}{1+x}\right\}_2 = 0.\]

Then
\[ \Delta = \{y\}_2\otimes x - \left\{ \frac{1-x}{1+x}\right\}_2  \otimes y + \left\{ \frac{x-1}{1+x}\right\}_2 \otimes y .\]

Let us write $ y = t^2$ so that
\begin{equation} \label{cond}
\frac{1+x}{1-x} = \pm \frac{t}{1-t^2}.
\end{equation}
Let us take the ``+"-case. Then
\[ x = \frac{t^2 + t -1}{-t^2+t+1}.\]

Thus we may write,
\[\Delta =\{t^2\}_2\otimes \frac{t^2 + t -1}{-t^2+t+1} -\left\{ \frac{1-t^2}{t} \right\}_2  \otimes t^2 + \left\{ \frac{t^2-1}{t}\right\}_2 \otimes t^2.\]
It is clear that a change in the sign of $t$ will simply change the sign of $\Delta$.

We will split the integration of $\Delta$ into two steps. In order to do that, write
\[ \Delta_1 = \{t^2\}_2\otimes \frac{t^2 + t -1}{-t^2+t+1}, \qquad \Delta_2 = 2\{t-t^{-1}\}_2\otimes t - 2\{t^{-1}-t\}_2\otimes t, \]
so that
\[ \Delta = \Delta_1 + \Delta_2.\]

We will first work with $\Delta_1$. Let $\varphi = \frac{1 + \sqrt{5}}{2}$, so $\varphi^2 - \varphi -1 = 0$.

By the five-term relation,
\begin{equation}\label{five1}
 \{\varphi t\}_2 + \{(\varphi-1)t\}_2 + \{1-t^2\}_2 + \left\{ \frac{1-\varphi t}{1-t^2} \right\}_2 + \left\{ \frac{1-(\varphi -1) t}{1-t^2} \right\}_2 =0,
\end{equation}
\begin{equation}\label{five2}
 \{-\varphi t\}_2 + \{(1-\varphi)t\}_2 + \{1-t^2\}_2 + \left\{ \frac{1+\varphi t}{1-t^2} \right\}_2 + \left\{ \frac{1+(\varphi-1) t}{1-t^2} \right\}_2 =0.
\end{equation}

Observe that we have 
\[ \Delta_1 = \{t^2\}_2 \otimes \frac{t^2 + t -1}{-t^2+t+1} =
\{t^2\}_2 \otimes \frac{(\varphi t-1)((\varphi-1)t+1)}{(\varphi t +1)((\varphi -1)t -1)}\]
\[ =\{t^2\}_2 \otimes \frac{1-\varphi t}{1- (\varphi -1)t } - \{t^2\}_2 \otimes \frac{1+\varphi t}{1+ (\varphi -1)t }.\]
Now let us apply the five-term relations (\ref{five1}) and (\ref{five2}):
\[ =
\{\varphi t\}_2 \otimes \frac{1-\varphi t}{1- (\varphi -1)t }+ \{(\varphi-1)t\}_2 \otimes \frac{1-\varphi t}{1- (\varphi -1)t } \]
\[+ \left\{ \frac{1-\varphi t}{1-t^2} \right\}_2 \otimes \frac{1-\varphi t}{1- (\varphi -1)t }+ \left\{ \frac{1-(\varphi -1) t}{1-t^2} \right\}_2 \otimes \frac{1-\varphi t}{1- (\varphi -1)t } \]
\[ - \{-\varphi t\}_2 \otimes \frac{1+\varphi t}{1+ (\varphi -1)t }- \{(1-\varphi)t\}_2 \otimes \frac{1+\varphi t}{1+ (\varphi -1)t } \]
\[- \left\{ \frac{1+\varphi t}{1-t^2} \right\}_2 \otimes \frac{1+\varphi t}{1+ (\varphi -1)t }- \left\{ \frac{1+(\varphi-1) t}{1-t^2} \right\}_2 \otimes \frac{1+\varphi t}{1+ (\varphi -1)t }.\]

Then we obtain
\[ \Delta_1 =  -\{1-\varphi t\}_2 \otimes (1 - \varphi t) +\{1-(\varphi-1) t\}_2 \otimes (1- (\varphi -1)t ) \]
\[+  \{1+\varphi t\}_2 \otimes (1 + \varphi t) -\{1+(\varphi-1) t\}_2 \otimes (1+ (\varphi -1)t ) \]
\[ - \{\varphi t\}_2 \otimes (1-(\varphi -1)t) +  \{-\varphi t\}_2 \otimes (1+(\varphi -1)t) +\{(\varphi -1)t\}_2\otimes(1-\varphi t) - \{(1-\varphi )t\}_2\otimes(1+\varphi t) \]
\[+ \left\{ \frac{1-\varphi t}{1-t^2} \right\}_2 \otimes \frac{1-\varphi t}{1- (\varphi -1)t }+ \left\{ \frac{1-(\varphi -1) t}{1-t^2} \right\}_2 \otimes \frac{1-\varphi t}{1- (\varphi -1)t } \]
\[- \left\{ \frac{1+\varphi t}{1-t^2} \right\}_2 \otimes \frac{1+\varphi t}{1+ (\varphi -1)t }- \left\{ \frac{1+(\varphi-1) t}{1-t^2} \right\}_2 \otimes \frac{1+\varphi t}{1+ (\varphi -1)t }.\]

Now we will  work with $\Delta_2$. By the five-term relation,
\[ \{ \varphi + t^{-1} \}_2 + \{1- (\varphi -1) t\}_2 + \{t - t^{-1}\}_2 + \left\{ \frac{1+(\varphi -1)t}{1-t^2}\right\}_2 + \left\{ \frac{(\varphi -1) t }{t - t^{-1}}\right \}_2 =0,\]

\[ \{ \varphi - t^{-1} \}_2 + \{1- (1-\varphi ) t\}_2 + \{t^{-1} - t\}_2 + \left\{ \frac{1-(\varphi -1)t}{1-t^2}\right\}_2 + \left\{ \frac{(\varphi -1) t }{t - t^{-1}}\right \}_2 =0,\]

\[ \{ 1 + \varphi t \}_2 + \{t^{-1} - (\varphi -1) \}_2 + \{t - t^{-1}\}_2 + \left\{ \frac{\varphi t}{t^{-1}-t}\right\}_2 + \left\{ \frac{1-\varphi t} {1 - t^2}\right \}_2 =0,\]

\[ \{ 1 - \varphi t \}_2 + \{-t^{-1} - (\varphi -1) \}_2 + \{t^{-1} - t\}_2 + \left\{ \frac{\varphi t}{t^{-1}-t}\right\}_2 + \left\{ \frac{1+\varphi t} {1 - t^2}\right \}_2 =0.\]

Applying the above equalities, we obtain
\[ \Delta_2 =  2\left\{ t - t^{-1} \right\}_2  \otimes t - 2\left\{ t^{-1} - t\right\}_2 \otimes t \]
\[ = \{ \varphi - t^{-1} \}_2 \otimes t - \{(1-\varphi ) t\}_2  \otimes t + \left\{ \frac{1-(\varphi -1)t}{1-t^2}\right\}_2  \otimes t\]
\[ - \{ \varphi + t^{-1} \}_2  \otimes t+ \{(\varphi -1) t\}_2  \otimes t - \left\{ \frac{1+(\varphi -1)t}{1-t^2}\right\}_2  \otimes t\]
\[ -\{  \varphi t \}_2  \otimes t- \{\varphi +t^{-1}  \}_2  \otimes t +   \left\{ \frac{1+\varphi t} {1 - t^2}\right \}_2 \otimes t\]
\[ +\{ - \varphi t \}_2  \otimes t+ \{\varphi - t^{-1}  \}_2  \otimes t- \left\{ \frac{1-\varphi t} {1 - t^2}\right \}_2 \otimes t.\]

We will now use the fact that we will integrate in a set where $|t|=1$.
Under those circumstances we have the following two identities (at the level of the differential form $\omega$):
 \begin{equation}
 \{(\varphi -1) t\}_2 \otimes t  =  \{(\varphi-1)\bar{t}\}_2 \otimes \bar{t} = \{\varphi t \}_2 \otimes t,  
 \end{equation}
 \begin{equation}
  \{(1-\varphi) t\}_2 \otimes t  =  \{-\varphi t \}_2 \otimes t. 
  \end{equation} 
Both identities depend on the fact that $|t|=1$ since we are conjugating and using that $\bar{t}=t^{-1}$.
We will work henceforth on the level of $\omega$ (i.e., in the coimage of $\eta_3(2)$). Thus,
\[ \Delta_2 = 2 \{ \varphi - t^{-1} \}_2 \otimes t - 2  \{ \varphi + t^{-1} \}_2  \otimes t\]
\[ + \left\{ \frac{1-(\varphi -1)t}{1-t^2}\right\}_2  \otimes t - \left\{ \frac{1+(\varphi -1)t}{1-t^2}\right\}_2  \otimes t +   \left\{ \frac{1+\varphi t} {1 - t^2}\right \}_2 \otimes t - \left\{ \frac{1-\varphi t} {1 - t^2}\right \}_2 \otimes t.\]

We will add $\Delta_1$ and $\Delta_2$. But first, let us note
\[\left\{ \frac{1-\varphi t}{1-t^2} \right\}_2 \otimes \frac{1-\varphi t}{t(1- (\varphi -1)t) } = \left\{ \frac{t^{-1}-\varphi }{t^{-1}-t} \right\}_2 \otimes \frac{t^{-1}-\varphi }{t-\varphi} - \left\{ \frac{1-\varphi t}{1-t^2} \right\}_2 \otimes(1-\varphi) \]
\[=-\left\{ \frac{t^{-1}-t} {t^{-1}-\varphi }\right\}_2 \otimes \frac{t^{-1}-\varphi }{t-\varphi} - \left\{ \frac{1-\varphi t}{1-t^2} \right\}_2 \otimes(1-\varphi)  \]
\[=- \left\{ \frac{t - \varphi} {t^{-1}-\varphi }\right\}_2 \otimes \frac{t-\varphi} {t^{-1}-\varphi }- \left\{ \frac{1-\varphi t}{1-t^2} \right\}_2 \otimes(\varphi-1).\]

\bigskip

Using this and similar identities, we obtain
\[ \Delta =-\{1-\varphi t\}_2 \otimes (1 - \varphi t) +\{1-(\varphi-1) t\}_2 \otimes (1- (\varphi -1)t ) \]
\[+  \{1+\varphi t\}_2 \otimes (1 + \varphi t) -\{1+(\varphi-1) t\}_2 \otimes (1+ (\varphi -1)t ) \]
\[ - \{\varphi t\}_2 \otimes (1-(\varphi -1)t) +  \{-\varphi t\}_2 \otimes (1+(\varphi -1)t) +\{(\varphi -1)t\}_2\otimes(1-\varphi t) - \{(1-\varphi )t\}_2\otimes(1+\varphi t) \]
\[-\left\{ \frac{t-\varphi }{t^{-1}-\varphi} \right\}_2 \otimes \frac{t-\varphi }{t^{-1}-\varphi} + \left\{ \frac{t-(\varphi -1) }{t^{-1}-(\varphi -1) } \right\}_2 \otimes \frac{t-(\varphi -1) }{t^{-1}-(\varphi -1) } \]
\[+ \left\{ \frac{t+\varphi }{t^{-1}+\varphi } \right\}_2 \otimes \frac{t+\varphi }{t^{-1}+\varphi } - \left\{ \frac{t+(\varphi -1) }{t^{-1}+(\varphi -1) } \right\}_2 \otimes \frac{t+(\varphi -1) }{t^{-1}+(\varphi -1) } \]
\[-\left\{ \frac{1-\varphi t}{1-t^2} \right\}_2 \otimes (\varphi -1) +  \left\{ \frac{1-(\varphi -1) t}{1-t^2} \right\}_2 \otimes \varphi \]
\[ +\left\{ \frac{1+\varphi t}{1-t^2} \right\}_2 \otimes ( \varphi -1) - \left\{ \frac{1+(\varphi -1)t}{1-t^2}\right\}_2  \otimes \varphi\]
\[+ 2 \{ \varphi - t^{-1} \}_2 \otimes t - 2  \{ \varphi + t^{-1} \}_2  \otimes t.\]

Now observe that
\[\left\{ \frac{1-\varphi t}{1-t^2} \right\}_2 \otimes \varphi  +  \left\{ \frac{1-(\varphi -1) t}{1-t^2} \right\}_2 \otimes \varphi  -\left\{ \frac{1+\varphi t}{1-t^2} \right\}_2 \otimes \varphi  - \left\{ \frac{1+(\varphi -1)t}{1-t^2}\right\}_2  \otimes \varphi\]
\[ =\{- \varphi t \}_2 \otimes \varphi  + \{(1- \varphi )t \}_2 \otimes \varphi -  \{\varphi t\}_2 \otimes \varphi - \{(\varphi-1) t\}_2 \otimes \varphi \]
by five-term relations (\ref{five1}) and (\ref{five2}).

Therefore (using $\frac{1}{\varphi-1} = \varphi$)
\[ \Delta =-\{1-\varphi t\}_2 \otimes (1 - \varphi t) +\{1-(\varphi-1) t\}_2 \otimes (1- (\varphi -1)t ) \]
\[+  \{1+\varphi t\}_2 \otimes (1 + \varphi t) -\{1+(\varphi-1) t\}_2 \otimes (1+ (\varphi -1)t ) \]
\[ - \{\varphi t\}_2 \otimes (1-(\varphi -1)t) +  \{-\varphi t\}_2 \otimes (1+(\varphi -1)t) +\{(\varphi -1)t\}_2\otimes(1-\varphi t) - \{(1-\varphi )t\}_2\otimes(1+\varphi t) \]
\[-\left\{ \frac{t-\varphi }{t^{-1}-\varphi} \right\}_2 \otimes \frac{t-\varphi }{t^{-1}-\varphi} + \left\{ \frac{t-(\varphi -1) }{t^{-1}-(\varphi -1) } \right\}_2 \otimes \frac{t-(\varphi -1) }{t^{-1}-(\varphi -1) } \]
\[+ \left\{ \frac{t+\varphi }{t^{-1}+\varphi } \right\}_2 \otimes \frac{t+\varphi }{t^{-1}+\varphi } - \left\{ \frac{t+(\varphi -1) }{t^{-1}+(\varphi -1) } \right\}_2 \otimes \frac{t+(\varphi -1) }{t^{-1}+(\varphi -1) } \]
\[+ \{- \varphi t \}_2 \otimes \varphi  + \{(1- \varphi )t \}_2 \otimes \varphi -  \{\varphi t\}_2 \otimes \varphi - \{(\varphi-1) t\}_2 \otimes \varphi \]
\[+ 2 \{ \varphi - t^{-1} \}_2 \otimes t - 2  \{ \varphi + t^{-1} \}_2  \otimes t.\]

Next we gather some similar terms together,
\[ \Delta =- \{1-\varphi t\}_2 \otimes (1 - \varphi t) +
\{1-(\varphi-1) t\}_2 \otimes (1- (\varphi -1)t ) \]
\[+  \{1+\varphi t\}_2 \otimes (1 + \varphi t) -
\{1+(\varphi-1) t\}_2 \otimes (1+ (\varphi -1)t ) \]
\[-\left\{ \frac{t-\varphi }{t^{-1}-\varphi} \right\}_2 \otimes \frac{t-\varphi }{t^{-1}-\varphi} + \left\{ \frac{t-(\varphi -1) }{t^{-1}-(\varphi -1) } \right\}_2 \otimes \frac{t-(\varphi -1) }{t^{-1}-(\varphi -1) } \]
\[+ \left\{ \frac{t+\varphi }{t^{-1}+\varphi } \right\}_2 \otimes \frac{t+\varphi }{t^{-1}+\varphi } - \left\{ \frac{t+(\varphi -1) }{t^{-1}+(\varphi -1) } \right\}_2 \otimes \frac{t+(\varphi -1) }{t^{-1}+(\varphi -1) } \]
\[ - \{\varphi t\}_2 \otimes (\varphi-t) +  \{-\varphi t\}_2 \otimes (\varphi +t) +\{(\varphi -1)t\}_2\otimes((\varphi-1)- t) - \{(1-\varphi )t\}_2\otimes((\varphi -1) + t) \]
\[+ 2 \{ \varphi - t^{-1} \}_2 \otimes t - 2  \{ \varphi + t^{-1} \}_2  \otimes t.\]

Observe that
\[ \{(\varphi-1) t \}_2 \otimes ((\varphi -1) - t ) =  \{(\varphi-1) t \}_2 \otimes ((\varphi -1) t ) + \{(\varphi-1) t \}_2 \otimes (\varphi  - t^{-1}  ) .\]
Now conjugate the elements of the second term (using $\varphi \in \Rset$)
\[=  \{(\varphi-1) t \}_2 \otimes ((\varphi -1) t ) + \{(\varphi-1) t^{-1} \}_2 \otimes (\varphi  - t  ) =  \{(\varphi-1) t \}_2 \otimes ((\varphi -1) t ) - \{\varphi t \}_2 \otimes (\varphi  - t  ).\]

Hence the last two lines of $\Delta$ above equal
\[ - \{\varphi t\}_2 \otimes (\varphi-t) +  \{-\varphi t\}_2 \otimes (\varphi +t) +\{(\varphi -1)t\}_2\otimes((\varphi-1)- t) - \{(1-\varphi )t\}_2\otimes((\varphi -1) + t) \]
\[+ 2 \{ \varphi - t^{-1} \}_2 \otimes t - 2  \{ \varphi + t^{-1} \}_2  \otimes t\]
\[ = - 2\{\varphi t\}_2 \otimes (\varphi-t) + 2 \{ \varphi - t^{-1} \}_2 \otimes t +  2\{-\varphi t\}_2 \otimes (\varphi +t) - 2  \{ \varphi + t^{-1} \}_2  \otimes t\]
\[+ \{(\varphi-1) t \}_2 \otimes ((\varphi -1) t ) -  \{(1- \varphi) t \}_2 \otimes ((1- \varphi ) t ).\]

We want to simplify the term $- \{\varphi t\}_2 \otimes (\varphi-t) +  \{ \varphi - t^{-1} \}_2 \otimes t$. On the one hand,
\[ - \{\varphi t \}_2 \otimes (\varphi  - t  ) = - \{\varphi t \}_2 \otimes t  + \{1- \varphi t \}_2 \otimes (1- \varphi  t^{-1}  ) .\]

On the other hand,
\[  \{ \varphi - t^{-1} \}_2 \otimes t = \{ \varphi - t \}_2 \otimes t^{-1} =  - \{ \varphi - t \}_2 \otimes \left ( \varphi - t \right )  + \{ \varphi - t \}_2 \otimes ( 1- \varphi t^{-1}).\]

By the five term relation
\[ \{1-\varphi t\}_2 + \{1-\varphi\}_2 + \{\varphi - t\}_2 - \{1-\varphi + t^{-1}\}_2 + \{(\varphi - 1) t \}_2 = 0,\]
but $\{1-\varphi\}_2$ corresponds to zero in the differential since it is a constant real number and $D(\Rset)=0$. We then get
\[  \{1- \varphi t \}_2 \otimes (1- \varphi  t^{-1}  )   +  \{ \varphi - t \}_2 \otimes ( 1- \varphi t^{-1})  = \{1- \varphi + t^{-1}\}_2 \otimes (1 -\varphi t^{-1}) - \{(\varphi-1)t\}_2  \otimes (1 -\varphi t^{-1})\]
\[=
 \{1- \varphi  + t^{-1}  \}_2   \otimes (1- \varphi  + t^{-1} ) -  \{1- \varphi  + t^{-1}  \}_2   \otimes (1-\varphi) \]
\[+ \{1-(\varphi - 1) t \}_2 \otimes (1-(\varphi - 1) t ) + \{(\varphi - 1) t \}_2 \otimes ((\varphi -1)t) .\]

Then
\[- \{\varphi t\}_2 \otimes (\varphi-t) +  \{ \varphi - t^{-1} \}_2 \otimes t = - \{\varphi t\}_2 \otimes t - \{\varphi - t\}_2 \otimes (\varphi - t) + \{1- \varphi  + t^{-1}  \}_2   \otimes (1- \varphi  + t^{-1} ) \]
\[-  \{1- \varphi  + t^{-1}  \}_2   \otimes (1-\varphi) + \{1-(\varphi - 1) t \}_2 \otimes (1-(\varphi - 1) t ) + \{(\varphi - 1) t \}_2 \otimes ((\varphi -1)t).\]

Analogously,
 \[ \{-\varphi t\}_2 \otimes (\varphi+t) -  \{ \varphi + t^{-1} \}_2 \otimes t =  \{-\varphi t\}_2 \otimes t + \{\varphi + t\}_2 \otimes (\varphi + t) - \{1- \varphi  - t^{-1}  \}_2   \otimes (1- \varphi  - t^{-1} ) \]
\[  +\{1- \varphi  - t^{-1}  \}_2   \otimes (1-\varphi) - \{1+(\varphi - 1) t \}_2 \otimes (1+(\varphi - 1) t ) - \{(1-\varphi ) t \}_2 \otimes ((1-\varphi)t).\]

Thus the last two lines of $\Delta$ above reduce to 
\[ - 2\{\varphi t\}_2 \otimes (\varphi-t) + 2 \{ \varphi - t^{-1} \}_2 \otimes t +  2\{-\varphi t\}_2 \otimes (\varphi +t) - 2  \{ \varphi + t^{-1} \}_2  \otimes t\]
\[+ \{(\varphi-1) t \}_2 \otimes ((\varphi -1) t ) -  \{(1- \varphi) t \}_2 \otimes ((1- \varphi ) t )  \]
\[ = - 2 \{\varphi t\}_2 \otimes \varphi t + 2\{\varphi t\}_2 \otimes \varphi -2 \{\varphi - t\}_2 \otimes (\varphi - t) + 2\{1- \varphi  + t^{-1}  \}_2   \otimes (1- \varphi  + t^{-1} ) \]
\[- 2 \{1- \varphi  + t^{-1}  \}_2   \otimes (1-\varphi) +2 \{1-(\varphi - 1) t \}_2 \otimes (1-(\varphi - 1) t ) +3\{(\varphi - 1) t \}_2 \otimes ((\varphi -1)t)  \]
\[ +  2\{-\varphi t\}_2 \otimes (-\varphi t) - 2\{-\varphi t\}_2 \otimes \varphi  + 2\{\varphi + t\}_2 \otimes (\varphi + t) -2 \{1- \varphi  - t^{-1}  \}_2   \otimes (1- \varphi  - t^{-1} ) \]
\[  +2\{1- \varphi  - t^{-1}  \}_2   \otimes (1-\varphi) -2 \{1+(\varphi - 1) t \}_2 \otimes (1+(\varphi - 1) t ) - 3\{(1-\varphi ) t \}_2 \otimes ((1-\varphi)t).\]

Next we will see that 
\[ \{ \varphi + t^{-1} \}_2 \otimes \varphi - \{\varphi -t^{-1} \}_2 \otimes \varphi + \{\varphi t\}_2 \otimes \varphi - \{-\varphi t\}_2 \otimes  \varphi   \]
corresponds to zero in the level of the differential form.

Using that $|t|=1$ and $\dd \arg \varphi =0$, the differential is
\[\frac{3 \omega}{\log \varphi} = \log| 1 - \varphi - t| \dd \log|\varphi + t| - \log | \varphi + t| \dd \log | 1 - \varphi - t| \]
\[ - \log| 1 - \varphi + t| \dd \log|\varphi - t| + \log | \varphi - t| \dd \log | 1 - \varphi + t| \]
\[ - \log \varphi \dd \log|1 - \varphi t| + \log \varphi \dd \log|1 + \varphi t| \] 

\[ =  \log| 1 + \varphi t| \dd \log|\varphi + t| -  \log \varphi \dd \log|\varphi + t| - \log | \varphi + t| \dd \log | 1 + \varphi  t| \]
\[ - \log| -1 + \varphi  t| \dd \log|\varphi - t| + \log \varphi \dd \log|\varphi - t| + \log | \varphi - t| \dd \log | 1 - \varphi  t| \]
\[ - \log \varphi \dd \log|t^{-1} - \varphi | + \log \varphi \dd \log|t^{-1} + \varphi | \] 

\[ =  \log| 1 + \varphi t| \dd \log|\varphi t^{-1}+ 1| - \log | \varphi t^{-1} + 1| \dd \log | 1 + \varphi  t| \]
\[ - \log| -1 + \varphi  t| \dd \log|\varphi t^{-1}- 1| + \log | \varphi t^{-1}-  1| \dd \log | 1 - \varphi  t| \]

\[=0\]

Finally the primitive of $\Delta$ is 
\[\Omega=  -\{1-\varphi t\}_3 + \{1-(\varphi-1) t\}_3 +  \{1+\varphi t\}_3 - \{1+(\varphi-1) t\}_3 \]
\[-\left\{ \frac{t-\varphi }{t^{-1}-\varphi} \right\}_3 + \left\{ \frac{t-(\varphi -1) }{t^{-1}-(\varphi -1) } \right\}_3 + \left\{ \frac{t+\varphi }{t^{-1}+\varphi } \right\}_3- \left\{ \frac{t+(\varphi -1) }{t^{-1}+(\varphi -1) } \right\}_3 \]
\[ - 2\{\varphi t\}_3 -2 \{\varphi - t\}_3 + 2\{1- \varphi  + t^{-1}  \}_3   + 2\{1-(\varphi - 1) t \}_3  +3\{(\varphi - 1) t \}_3 \]
\[ +2  \{-\varphi t\}_3+ 2\{\varphi + t\}_3 -2 \{1- \varphi  - t^{-1}  \}_3   -2 \{1+(\varphi - 1) t \}_3 - 3\{(1-\varphi ) t \}_3,\]
which is
\[ =  -\{1-\varphi t\}_3  +  \{1+\varphi t\}_3  \]
\[-\left\{ \frac{t-\varphi }{t^{-1}-\varphi} \right\}_3 + \left\{ \frac{t-(\varphi -1) }{t^{-1}-(\varphi -1) } \right\}_3 + \left\{ \frac{t+\varphi }{t^{-1}+\varphi } \right\}_3- \left\{ \frac{t+(\varphi -1) }{t^{-1}+(\varphi -1) } \right\}_3 \]
\[ - 2\{\varphi t\}_3 -2 \{\varphi - t\}_3 + 2\{1- \varphi  + t^{-1}  \}_3   + 3\{1-(\varphi - 1) t \}_3  +3\{(\varphi - 1) t \}_3 \]
\[ +2  \{-\varphi t\}_3+ 2\{\varphi + t\}_3 -2 \{1- \varphi  - t^{-1}  \}_3   -3 \{1+(\varphi - 1) t \}_3 - 3\{(1-\varphi ) t \}_3.\]

Let us use that 
\[\{x\}_3 +  \{ 1 - x\}_3 + \left \{ 1 - \frac{1}{x} \right \}_3 = \{1\}_3 \]
and the fact that $\{x\}_3=\{\bar{x}\}_3$ at the level of the differential.

We obtain
 \[\Omega= - 4\{1-\varphi t\}_3  + 4 \{1+\varphi t\}_3  \]
\[-\left\{ \frac{t-\varphi }{t^{-1}-\varphi} \right\}_3 + \left\{ \frac{t-(\varphi -1) }{t^{-1}-(\varphi -1) } \right\}_3 + \left\{ \frac{t+\varphi }{t^{-1}+\varphi } \right\}_3- \left\{ \frac{t+(\varphi -1) }{t^{-1}+(\varphi -1) } \right\}_3 \]
\[ - 2\{\varphi t\}_3 -2 \{\varphi - t\}_3 + 2\{1- \varphi  + t^{-1}  \}_3   +2  \{-\varphi t\}_3+ 2\{\varphi + t\}_3 -2 \{1- \varphi  - t^{-1}  \}_3 .\]

Let us note that the poles of $z$ occur with $x=1$, which easily implies $\Delta=0$. Analogously, $\Delta =0$ for $y=1$ or $x=-1$ which correspond to the zeros of $z$.  

We need to describe the integration path. If we let $x=\e^{2\ii\alpha}$, with $-\frac{\pi}{2} \leq \alpha \leq \frac{\pi}{2}$, and $t=\e^{\ii \beta}$, with  $-\frac{\pi}{2}\leq \beta \leq \frac{\pi}{2}$. Then condition (\ref{cond}) translates into
\[ \tan \alpha = \pm 2 \sin \beta.\]

\begin{figure} 
\centering
\includegraphics[width=12pc]{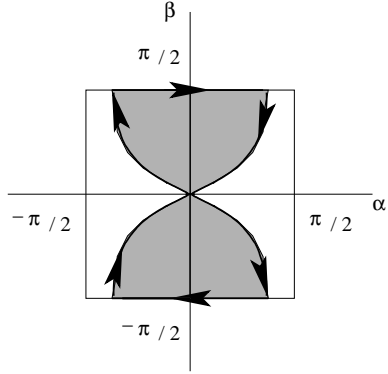}
\caption{\label{mat3condon}Integration set for $(1-y)(1+x)-z(1-x)$.}
\end{figure}

The boundaries of the above condition are met when $\sin \beta = \pm 1$, corresponding to $t=\pm \ii$. Also, by inspection of Eq. (\ref{cond}) we see that $t=\pm 1$ and $x=1$ are also critical points.  After carefully analyzing the situation (see Fig. \ref{mat3condon}), and taking into account that $\Delta =0$ whenever $x=1$ or $t=\pm 1$, we conclude that we need to integrate $\Delta$ with $-\frac{\pi}{2}\leq \beta \leq 0$  and $\frac{\pi}{2}\geq \beta \geq 0$. Since $\Omega$ is 0 when evaluated for $t=\pm \ii$, and $\Omega(1) = - \Omega(-1)$, we obtain
\[ 4 \pi^2 m(P) = 4 \Omega\,(1)\]
\[=4( -4 \mathcal{L}_3(1 - \varphi) + 4 \mathcal{L}_3(1+\varphi) - 2\mathcal{L}_3(\varphi)+ 2 \mathcal{L}_3(- \varphi) - 2\mathcal{L}_3(\varphi-1) + 2\mathcal{L}_3(\varphi +1)\]
\[ + 2 \mathcal{L}_3( 2 - \varphi) - 2 \mathcal{L}_3(-\varphi) ) )\] 
\[ = 4( -4 \mathcal{L}_3(1 - \varphi) + 6 \mathcal{L}_3(1+\varphi) - 2\mathcal{L}_3(\varphi) - 2\mathcal{L}_3(\varphi-1)  + 2 \mathcal{L}_3( 2 - \varphi)  ). \]

Now we use that
\[ \{1-\varphi\}_3 = \{-\varphi^{-1}\}_3=\{-\varphi\}_3, \qquad \{\varphi - 1\}_3 = \{\varphi\}_3,\]
in order to obtain
\[\pi^2 m(P) = 6 \mathcal{L}_3(1+ \varphi) - 4 \mathcal{L}_3(\varphi  ) - 4 \mathcal{L}_3(-\varphi) + 2 \mathcal{L}_3(2-\varphi).\]

Then we use
\[ \{ \varphi -1 \}_3 + \{2-\varphi\}_3 + \left\{ 1 - \frac{1}{\varphi-1}\right\}_3 = \{1\}_3, \]
\[ \{ \varphi  \}_3 + \{2-\varphi\}_3 + \left\{ - \varphi \right\}_3 = \{1\}_3,\]
and 
\[ \{- \varphi  \}_3 + \{1+\varphi\}_3 + \left\{ 1 + \frac{1}{\varphi}\right\}_3 = \{1\}_3, \]
\[ \{- \varphi  \}_3 + \{1+\varphi\}_3 + \left\{ \varphi\right\}_3 = \{1\}_3.\]

Thus we obtain
\[ \pi^2 m(P)  =  8 \mathcal{L}_3(1) - 12\mathcal{L}_3(\varphi) - 12 \mathcal{L}_3(-\varphi).\]

But 
\[ 4\{\varphi\}_3 + 4 \{-\varphi\}_3 = \{\varphi^2\}_3 = \{\varphi + 1\}_3 = \{1\}_3 - \{\varphi\}_3 - \{-\varphi\}_3\]
implies that
\[ \{\varphi\}_3 + \{-\varphi\}_3 = \frac{1}{5} \{1\}_3.\]

Finally we recover Condon's result
\begin{equation}
 m(P) = \frac{28}{5 \pi^2} \zeta(3).
\end{equation}

Let us note that the computation involves $\sqrt{5}$, although the final result does not. The same situation applies to Condon's computation in \cite{C} (in fact, he also needs to use $\varphi$), so it seems that this step is unavoidable as in the case of Smyth's example from Section \ref{Smyth2}.

Moreover, if we look at the rational number in each of the three-variable formulas we will notice that the denominators are always 1 or 2 except for these two examples of Smyth and Condon.
In Smyth's example the denominator is a 3 and in Condon's example the denominator is a 5. Both the appearance of 3 as a denominator together with $2\{3\}_3-\{-3\}_3$ in the first case, and 5 as a denominator together with  elements in $\mathcal{B}(\Qset(\sqrt{5}))$ seem to be related phenomena. Gangl suggested that a possible explanation may be that there are torsion elements  that we may be overlooking when we tensor by $\Qset$ (as the regulator maps kill torsion).

Lastly, this is the first example where we used non-rational functional equations for the polylogarithms, namely $D(\bar{x}) = -D(x)$ and $\mathcal{L}_3(\bar{x})=\mathcal{L}_3(x)$. This suggests that for us, the best definition for $\mathcal{R}_n(F)$ may be ``all the functional equations of $\mathcal{L}_n$" instead of ``all the rational functional equations of $\mathcal{L}_n$". Of course, the second definition has the advantage of being more concrete and should be easier to handle.

\subsection{An example in four variables}
We will study an example in four variables:
\[ z+ \frac{(1+x)(1+x_1)}{(1+y)(1-x_1)}.\]
This example was first computed in \cite{L} in terms of multiple polylogarithms. It was later observed \cite{L4} that this particular combination of multiple polylogarithms has a simpler expression in terms of a Dirichlet L-series. Here we relate this Mahler measure directly to the Dirichlet L-series.

The wedge product is
\[ x\wedge y\wedge x_1 \wedge z = x\wedge y\wedge x_1 \wedge \frac{(1+x)(1+x_1)}{(1+y)(1-x_1)}\]
\[ = x\wedge (1+x) \wedge y \wedge x_1 +y\wedge (1+y) \wedge x\wedge x_1 +x_1 \wedge (1+x_1) \wedge x \wedge y  - x_1 \wedge (1-x_1) \wedge x \wedge y. \]
We integrate and obtain
\[\Delta= \{-x\}_2 \otimes y \wedge x_1 +\{-y\}_2 \otimes x\wedge x_1 +\{-x_1\}_2\otimes x\wedge y -\{x_1\}_2\otimes x \wedge y.\]

It  is more convenient to eliminate the variable $x_1$ than $z$. Notice that
\[x_1=\frac{-(1+x)-z(1+y)}{(1+x)-z(1+y)}.\]
Then the condition whose intersection with $\TT^4$ describes $|\partial \Gamma|$ becomes
\[  (1+x)(1+y)(x+yz^2) =0.\]

When $x=-1$ or $y=-1$, $\eta_4(3) (\Delta) = 0$. Now suppose $x=-yz^2$.
Then we need to integrate
\[\Delta =  \{yz^2\}_2 \otimes y \wedge x_1 +\{-y\}_2 \otimes yz^2 \wedge x_1 +\{-x_1\}_2\otimes z^2\wedge y -\{x_1\}_2\otimes z^2 \wedge y.\]
Consider the first two terms,
\[ \Delta_1= \{yz^2\}_2 \otimes y \wedge x_1 +\{-y\}_2 \otimes yz^2 \wedge x_1 \]
\[= \{yz^2\}_2 \otimes yz^2 \wedge x_1 - 2 \{yz^2\}_2 \otimes z \wedge x_1 + \{-y\}_2 \otimes y\wedge x_1 + 2\{-y\}_2 \otimes z\wedge x_1.\]

Use the five-term relation,
\[ \{yz^2\}_2 + \left\{ -\frac{1}{y}\right\}_2 + \{1+z^2\}_2 +\left\{ \frac{1-yz^2}{1+z^2}\right\}_2 + \left\{\frac{y+1}{y(1+z^2)}\right\}_2 =0 ,\]
together with the fact that
\[ \alpha:= \frac{x_1-1}{1+x_1} = \frac{1-yz^2}{z(1+y)}\]
then
\[\{-y\}_2 \otimes z\wedge x_1- \{yz^2\}_2 \otimes z \wedge x_1 \]
\[= -\{-z^2\}_2 \otimes z \wedge x_1 + \left\{ \frac{1-yz^2}{1+z^2}\right\}_2\otimes z\wedge x_1 +  \left\{\frac{y+1}{y(1+z^2)}\right\}_2\otimes z\wedge x_1\]
\[= -\{-z^2\}_2 \otimes z \wedge x_1 + \left\{ -\frac{z^2(1+y)}{1-yz^2}\right\}_2\otimes z\wedge x_1 +  \left\{\frac{1-yz^2}{1+y}\right\}_2\otimes z\wedge x_1\]
\[= -\{-z^2\}_2 \otimes z \wedge x_1 + \left\{ -\frac{z}{\alpha}\right\}_2\otimes z\wedge x_1 +  \left\{z\alpha\right\}_2\otimes z\wedge x_1.\]
Conjugate and invert the second term (noting that $\arg \alpha = \pm \frac{\pi}{2}$ since $|x_1|=1$)
\[= -\{-z^2\}_2 \otimes z \wedge x_1 + 2  \left\{z\alpha\right\}_2\otimes z\wedge x_1.\]

Consider the last two terms of $\Delta$:
\[ \Delta_2 = \{-x_1\}_2\otimes z^2\wedge y -\{x_1\}_2\otimes z^2 \wedge y .\]

Note that
\[ y = \frac{1-\alpha z}{z(z+\alpha)},\]
and recall the five term relation:
\[\{x_1\}_2 - \{-x_1\}_2 = \left\{\frac{x_1-1}{1+x_1}\right\}_2 - \left\{\frac{1-x_1}{1+x_1}\right\}_2 = \{\alpha\}_2-\{-\alpha\}_2.\]
Thus
\[ \Delta_2= \{-\alpha\}_2\otimes z^2\wedge  \frac{1-\alpha z}{z(z+\alpha)}-\{\alpha\}_2\otimes z^2 \wedge \frac{1-\alpha z}{z(z+\alpha)}.\]
Conjugating the first term, and inverting,
\[ \Delta_2= -2\{\alpha\}_2\otimes z^2 \wedge \frac{1-\alpha z}{z(z+\alpha)}. \]
Thus
\[\Delta= \{yz^2\}_2 \otimes yz^2 \wedge x_1+ \{-y\}_2 \otimes y\wedge x_1-\{-z^2\}_2 \otimes z^2 \wedge x_1\]
\[+4 \left\{z\alpha\right\}_2\otimes z\wedge \frac{1+\alpha}{1-\alpha} - 4 \{\alpha\}_2\otimes z \wedge  \frac{1-\alpha z}{z+\alpha}.\]

By inverting and dividing by $z$:
\[\{\alpha\}_2\otimes z \wedge  \frac{1-\alpha z}{z+\alpha}= -\{\alpha\}_2\otimes z \wedge  \frac{1+\alpha z^{-1} }{1-\alpha z}.\]

Now notice that
\[\left\{z\alpha\right\}_2\otimes z\wedge \frac{1+\alpha}{1-\alpha}+  \{\alpha\}_2\otimes z \wedge  \frac{1+\alpha z^{-1}}{1-\alpha z}\]
\[=\left\{z\alpha\right\}_2\otimes z\wedge \frac{1+\alpha}{1-z\alpha} -\left\{z\alpha\right\}_2\otimes z\wedge \frac{1-\alpha}{1-z\alpha} + \{\alpha\}_2\otimes z \wedge  \frac{1-\alpha}{1-\alpha z}- \{\alpha\}_2\otimes z \wedge  \frac{1-\alpha}{1+\alpha z^{-1}}.\]
Conjugate and invert $z$ in the last term
\[\left\{z\alpha\right\}_2\otimes z\wedge \frac{1+\alpha}{1-z\alpha} -\left\{z\alpha\right\}_2\otimes z\wedge \frac{1-\alpha}{1-z\alpha} + \{\alpha\}_2\otimes z \wedge  \frac{1-\alpha}{1-\alpha z}- \{-\alpha\}_2\otimes z \wedge  \frac{1+\alpha}{1-\alpha z}.\]

Now use
\[\{\alpha\}_2+\{z\}_2+\{1-\alpha z\}_2+ \left\{\frac{1-\alpha}{1-\alpha z}\right\}_2 +  \left\{\frac{1-z}{1-\alpha z}\right\}_2=0,\]

\[\{-\alpha\}_2+\{-z\}_2+\{1-\alpha z\}_2+ \left\{\frac{1+\alpha}{1-\alpha z}\right\}_2 +  \left\{\frac{1+z}{1-\alpha z}\right\}_2=0.\]

Then we obtain
\[-\{z\}_2\otimes z\wedge \frac{1-\alpha}{1-z\alpha} - \left\{\frac{1-\alpha}{1-\alpha z}\right\}_2
\otimes z\wedge \frac{1-\alpha}{1-z\alpha} + \left\{\frac{z(1-\alpha)}{1-\alpha z}\right\}_2\otimes z\wedge \frac{1-\alpha}{1-z\alpha} \]
\[+\{-z\}_2\otimes z\wedge \frac{1+\alpha}{1-z\alpha} + \left\{\frac{1+\alpha}{1-\alpha z}\right\}_2 \otimes z\wedge \frac{1+\alpha}{1-z\alpha} - \left\{-\frac{z(1+\alpha)}{1-\alpha z}\right\}_2\otimes z\wedge \frac{1+\alpha}{1-z\alpha}.\]

Then $\Delta$ integrates (on the component $\gamma_0 = \{ x=-yz^2\}$ of $|\partial \Gamma|$) to
\[\Omega= \{yz^2\}_3 \otimes x_1+ \{-y\}_3 \otimes x_1-\{-z^2\}_3 \otimes x_1 \]
\[- 4\{z\}_3\otimes \frac{1-\alpha}{1-z\alpha} +4 \left\{\frac{1-\alpha}{1-\alpha z}\right\}_3
\otimes z - 4 \left\{\frac{z(1-\alpha)}{1-\alpha z}\right\}_3\otimes z \]
\[+4\{-z\}_3\otimes \frac{1+\alpha}{1-z\alpha} -4 \left\{\frac{1+\alpha}{1-\alpha z}\right\}_3 \otimes z +4 \left\{-\frac{z(1+\alpha)}{1-\alpha z}\right\}_3\otimes z .\]

Now the boundary of $\gamma_0$ is $x=-1$ or $y=-1$. When $x=-1$ we obtain $y=-1$ and $z= \pm \ii$ or $\alpha=0$ and $x_1=1$. If $y=-1$,
\[ \Omega= \{1\}_3 \otimes x_1- 4\{z\}_3\otimes \frac{1-\alpha}{1-z\alpha} +4\{-z\}_3\otimes \frac{1+\alpha}{1-z\alpha}.\]
If we take into account that $\alpha$ moves in the imaginary axis when $x_1$ moves in the unit circle and that we are evaluating on $z=\pm \ii$, then we just need to evaluate $\eta_4(2)$ on
\[\Omega = \{1\}_3 \otimes x_1. \]

If $\alpha=0$,
\[\Omega=   -4 \left\{z\right\}_3\otimes z +4 \left\{-z\right\}_3\otimes z,\]
integrating to
\[  \Sigma_1= -4 \{z\}_4+4\{-z\}_4,\]
and the boundary now is $z^2=-1$.

When $y=-1$, then $z= \pm \ii$ and $x=-1$ or $x_1=-1$ and $\alpha=\infty$. In the first case  we obtain, as before
\[ \Omega= \{1\}_3 \otimes x_1.\]
In the second case we get,
\[\Omega= 8\{z\}_3\otimes z-8\{-z\}_3\otimes z \]
integrating to 
\[  \Sigma_2=  8 \{z\}_4-8\{-z\}_4.\]
Now we need to take into account the opposite orientations that will cancel the terms with $ \Omega= \{1\}_3 \otimes x_1$ and add up to
\[\Sigma =12\{z\}_4-12\{-z\}_4,\]
which must be evaluated between $\ii$ and $-\ii$. This yields the final result
\begin{equation}
m(P) = \frac{24}{\pi^3} \Lf(\chi_{-4}, 4).
\end{equation}

\section{Generalized Mahler measure}

In this section we will apply the algebraic integration to the computation of generalized Mahler measures. 

We follow the notation from Eq. (\ref{defn2}). Let us fix $P \in \Rset[x^{\pm 1}]$ and take $f_j = P(x_j)$ (so in particular $f_j \in \Rset[x_j^{\pm 1}] \subset \Rset[x_1^{\pm 1}, \dots, x_n^{\pm 1}] $). 

Suppose, moreover, that $|P|$ is a monotonic function for $0 \leq \arg x \leq \pi$. Write $x_j = \e^{2\pi \ii \theta_j}= \e(\theta_j)$. Then
\[m(P(x_1),\dots, P(x_n)) =\int_{-\frac{1}{2}}^{\frac{1}{2}}\dots \int_{-\frac{1}{2}}^{\frac{1}{2}} \max\{\log|P(\e(\theta_1))|, \dots, \log|P(\e(\theta_n))|\} \dd \theta_1 \dots \dd \theta_n\]
\[=2^n n! \int_{0\leq \theta_n \leq \dots \leq \theta_1 \leq \frac{1}{2}} \log|P(\e(\theta_1))| \dd \theta_1 \dots \dd \theta_n   \]
\[=\frac{n!}{(\pi \ii)^n} \int_ \Gamma\eta_{n+1}(n+1)(P(x_1),x_1,\dots,x_n)\]
where $\Gamma = \{0\leq \arg(x_n) \leq \dots \leq \arg(x_1) \leq \pi\}$.

We proceed to compute some examples.

\subsection{The case of $P=1-x$}
Let $P=1-x$, then $|P(\e(\theta))|= 2 \left |\sin \left(\pi\theta\right ) \right |$, so it is monotonic on $[0,\pi]$. This  case was computed by Gon and Oyanagi in \cite{GO} with the use of multiple sine functions.

We start by evaluating $\eta_{n+1}(n+1)$ on 
\[ (1-x_1) \wedge x_1 \wedge \dots \wedge x_n.\]
Now this corresponds to an exact form, which we integrate and obtain
\[-\{x_1\}_2 \otimes x_2 \wedge \dots \wedge x_n.\]
The boundary $|\partial \Gamma|$ consists of $x_1= -1$ and $x_1=x_2$ for our purposes (the remaining restrictions are $0$, e.g., $x_2=x_3$ implies $x_2\wedge x_3 =0$). Evaluating,
\[ -\{-1\}_2 \otimes x_2 \wedge \dots \wedge x_n + \{x_2\}_2 \otimes x_2 \wedge x_3 \wedge \dots \wedge x_n.\]
By integrating the term in the right and evaluating in the boundary (of $|\partial \Gamma| \cap \{x_1=x_2\}$) $x_2=-1$ and $x_2=x_3$, we obtain,
\[ -\{-1\}_2 \otimes x_2 \wedge \dots \wedge x_n \oplus\{-1\}_3 \otimes x_3 \wedge \dots \wedge x_n - \{x_3\}_3 \otimes x_3 \wedge \dots \wedge x_n.\]
Here, as well as in the following pages, $\oplus$ and $\bigoplus$ will denote formal versions of $+$ and $\sum$. We will use this notation whenever we have a combination of forms of different dimensions. Eventually, each of these forms will be integrated in a simplex of the right dimension, and the results of each of these integrals will be added to obtain the final result.

We continue with the same procedure, until we reach
\[ (-1)^{n+1} \{1\}_{n+1} \oplus\bigoplus_{k=2}^{n+1} (-1)^{k+1} \{-1\}_k \otimes x_k \wedge \dots \wedge x_n.\]
Now we need to integrate each term. First observe that since the $\mathcal{L}_k$ is trivial in the real numbers for $k$ even, we just need to consider
\[  -\{1\}_{n+1} \oplus \bigoplus_{k=3,\, k \,\mathrm{odd}}^{n+1}  \{-1\}_k \otimes x_k \wedge \dots \wedge x_n,\]
where it is understood that the sum is up to $n+1$ if $n$ is even or $n$ if $n$ is odd.

The sum yields
\[\sum_{k=3,\, k \,\mathrm{odd}}^{n+1} \frac{(\ii \pi)^{n-k+1}}{(n-k+1)!}\Li_k(-1) =  \sum_{k=3,\, k \,\mathrm{odd}}^{n+1} \frac{(\ii \pi)^{n-k+1}(1-2^{k-1})}{(n-k+1)!2^{k-1}} \zeta(k).\]
Finally, we get
\begin{eqnarray}
m(1-x_1,\dots,1-x_{2m}) &=& \frac{(-1)^{m+1}(2m)!}{\pi^{2m}} \zeta(2m+1) \nonumber \\
&&+(2m)!\sum_{j=1}^m (-1)^{j} \frac{(1-2^{2j})}{(2m-2j)!(2\pi)^{2j}} \zeta(2j+1),\\
\label{limit} m(1-x_1,\dots,1-x_{2m-1}) &=& (2m-1)! \sum_{j=1}^{m-1} (-1)^{j} \frac{(1-2^{2j})}{(2m-2j-1)!(2\pi)^{2j}} \zeta(2j+1),
\end{eqnarray}
recovering the result of \cite{GO}.
\subsection{The example with $P=\frac{1-x}{1+x}$}
Although $P=\frac{1-x}{1+x}$ is a rational function, we can still apply the method to this case.
Notice that $|P(\e(\theta))| = \left | \tan \left ( \pi\theta)\right ) \right|$ so it is monotonic. We start considering
\[ \frac{1-x_1}{1+x_1} \wedge x_1 \wedge \dots \wedge x_n. \]
We can find the primitive easily as 
\[\{-x_1\}_2 \otimes x_2 \wedge \dots \wedge x_n-\{x_1\}_2 \otimes x_2 \wedge \dots \wedge x_n.\]
The boundary consists of $x_1=-1$ and $x_1=x_2$, and thus we obtain
\[ \{1\}_2 \otimes x_2 \wedge \dots \wedge x_n-\{-1\}_2 \otimes x_2 \wedge \dots \wedge x_n\]
\[- \{-x_2\}_2 \otimes x_2 \wedge \dots \wedge x_n+\{x_2\}_2 \otimes x_2 \wedge \dots \wedge x_n.\]
If we continue the integration steps, we reach
\[(-1)^{n+1} (\{1\}_{n+1}-\{-1\}_{n+1}) \oplus \bigoplus_{k=2}^{n+1} (-1)^{k}(\{1\}_k \otimes x_k \wedge \dots \wedge x_n- \{-1\}_k \otimes x_k \wedge \dots \wedge x_n).\]
Taking into account the parity, we just need to consider
\[\{-1\}_{n+1}-\{1\}_{n+1} \oplus\bigoplus_{k=3, \,k\,\mathrm{odd}}^{n+1} \{-1\}_k \otimes x_k \wedge \dots \wedge x_n- \{1\}_k \otimes x_k \wedge \dots \wedge x_n.\]

The sum yields
\[ \sum_{k=3,\, k \,\mathrm{odd}}^{n+1} \frac{(\ii\pi)^{n-k+1}}{(n-k+1)!}(\Li_k(-1)-\Li_k(1)) =  \sum_{k=3,\, k \,\mathrm{odd}}^{n+1} \frac{(\ii \pi)^{n-k+1}(1-2^{k})}{(n-k+1)!2^{k-1}}\zeta(k).\]

Then we reach the result
\begin{eqnarray}
 m\left( \frac{1-x_1}{1+x_1}, \dots, \frac{1-x_{2m}}{1+x_{2m}}\right) &= &\frac{(-1)^{m}(2m)!(1-2^{2m+1})}{(2\pi)^{2m}} \zeta(2m+1) \nonumber \\
&&+(2m)!\sum_{j=1}^m (-1)^{j} \frac{(1-2^{2j+1})}{(2m-2j)!(2\pi)^{2j}} \zeta(2j+1),\\
m\left(\frac{1-x_1}{1+x_1},\dots,\frac{1-x_{2m-1}}{1+x_{2m-1}}\right) &=& (2m-1)! \sum_{j=1}^{m-1} (-1)^{j} \frac{(1-2^{2j+1})}{(2m-2j-1)!(2\pi)^{2j}} \zeta(2j+1).\nonumber\\
&&
\end{eqnarray}
We remark the specific case
\[m\left( \frac{1-x_1}{1+x_1}, \frac{1-x_{2}}{1+x_{2}}\right) = \frac{7}{\pi^2} \zeta(3).\]
It corresponds to the Mahler measure of  $ \frac{1-x_1}{1+x_1}+z \frac{1-x_{2}}{1+x_{2}}$. Thus we recover formula (\ref{eq:mine}).

\subsection{An example involving the golden ratio}
We will consider the case of $P=1+x-x^{-1}$. Observe that 
\[|P(\e(\theta))| = |1+\ii 2 \sin (2\pi\theta)| = \sqrt{1+4 \sin^2 (2\pi \theta)}.\]
 In this case we need to integrate with $0 \leq \arg x_i \leq \frac{\pi}{2}$ (since that is when $P\circ \e$ is monotonic). Then
\[m(P(x_1),\dots, P(x_n)) =n! \left( \frac{2}{\pi\ii}\right)^n \int_{0\leq \arg(x_n) \leq \dots \leq \arg(x_1) \leq \frac{\pi}{2}} \eta_{n+1}(n+1)(P(x_1),x_1,\dots,x_n).\]
We start with
\[ (1+x_1-x_1^{-1}) \wedge x_1 \wedge \dots \wedge x_n = (x_1^2+x_1-1) \wedge x_1 \wedge \dots \wedge x_n \]
\[ =(1-\varphi^{-1}x_1) \wedge x_1 \wedge \dots \wedge x_n + (1+\varphi x_1) \wedge x_1 \wedge \dots \wedge x_n, \]
where $\varphi = \frac{-1+\sqrt{5}}{2}$.

Integrating, we obtain
\[ -\{\varphi^{-1} x_1\}_2 \otimes x_2 \wedge \dots \wedge x_n - \{-\varphi x_1\}_2 \otimes x_2 \wedge \dots \wedge x_n  \oplus \left( \frac{\varphi -x_1}{1+\varphi x_1}\right) \wedge \varphi \wedge x_2 \wedge \dots \wedge x_n.\]
The boundary is given by $x_1=\ii$ and $x_1=x_2$, thus,
\[ -\{\varphi^{-1} \ii \}_2 \otimes x_2 \wedge \dots \wedge x_n - \{-\varphi \ii\}_2 \otimes x_2 \wedge \dots \wedge x_n +  \{\varphi^{-1} x_2\}_2 \otimes x_2 \wedge \dots \wedge x_n + \{-\varphi x_2\}_2 \otimes x_2 \wedge \dots \wedge x_n\]  
\[ \oplus \left( \frac{\varphi -x_1}{1+\varphi x_1}\right) \wedge \varphi \wedge x_2 \wedge \dots \wedge x_n.\]

Integrating once more, we obtain
\[ -\{\varphi^{-1} \ii \}_2 \otimes x_2 \wedge \dots \wedge x_n - \{-\varphi \ii\}_2 \otimes x_2 \wedge \dots \wedge x_n \oplus  \{\varphi^{-1} x_2\}_3 \otimes x_3 \wedge \dots \wedge x_n + \{-\varphi x_2\}_3 \otimes x_3 \wedge \dots \wedge x_n\] 
\[\oplus \{\varphi^{-1} x_2\}_2 \otimes  \varphi \wedge x_3 \wedge \dots \wedge x_n - \{-\varphi x_2\}_2 \otimes  \varphi \wedge x_3 \wedge \dots \wedge x_n\oplus \left( \frac{\varphi -x_1}{1+\varphi x_1}\right) \wedge \varphi \wedge x_2 \wedge \dots \wedge x_n.\]

Continuing with the same procedure, we reach
\[(-1)^{n+1} (\{\varphi^{-1}\}_{n+1} + \{-\varphi\}_{n+1}) \oplus\bigoplus_{k=2}^{n+1} (-1)^{k+1}(\{\varphi^{-1} \ii \}_k \otimes x_k \wedge \dots \wedge x_n + \{-\varphi \ii\}_k \otimes x_k \wedge \dots \wedge x_n)\] 
\[\oplus \left( \frac{\varphi -x_1}{1+\varphi x_1}\right) \wedge \varphi \wedge x_2 \wedge \dots \wedge x_n\]
\[\oplus \bigoplus_{k=2}^{n} (-1)^k (\{\varphi^{-1} x_k\}_k \otimes  \varphi \wedge x_{k+1} \wedge \dots \wedge x_n - \{-\varphi x_k\}_k \otimes  \varphi \wedge x_{k+1} \wedge \dots \wedge x_n).\]

Taking into account properties of $\mathcal{L}_k$ according to the parity of $k$,
\[=- \{\varphi\}_{n+1} - \{-\varphi\}_{n+1} \oplus \bigoplus_{k=3,\,k\,\mathrm{odd}}^{n+1} 2 \{-\varphi \ii\}_k \otimes x_k \wedge \dots \wedge x_n\] 
\[\oplus \left( \frac{\varphi -x_1}{1+\varphi x_1}\right) \wedge \varphi \wedge x_2 
\wedge \dots \wedge x_n\]
\begin{equation}\label{eq:golden}
\oplus \bigoplus_{k=2}^{n} (-1)^k (\{\varphi^{-1} x_k\}_k \otimes  \varphi \wedge x_{k+1} \wedge \dots \wedge x_n - \{-\varphi x_k\}_k \otimes  \varphi \wedge x_{k+1} \wedge \dots \wedge x_n).
\end{equation}
The first sum in (\ref{eq:golden}) yields
\[2\sum_{k=3,\, k \,\mathrm{odd}}^{n+1} \frac{(\ii \pi)^{n-k+1}}{2^{n-k+1}(n-k+1)!}\mathcal{L}_k(\varphi \ii)= \sum_{k=3,\, k \,\mathrm{odd}}^{n+1} \frac{(\ii \pi)^{n-k+1}}{2^{n}(n-k+1)!}\mathcal{L}_k(-\varphi^2).\]
The second line in (\ref{eq:golden}) yields
\[ -(-1)^n \log \varphi \left( \int_\ii^1 \ii\,{\im}\left( \frac{\dd x_1}{x_1-\varphi} \right) \circ \frac{\dd x_2}{x_2} \circ \dots \circ \frac{\dd x_n}{x_n}  - \int_\ii^1\ii\,{\im}\left( \frac{\dd x_1}{x_1+ \varphi^{-1}} \right) \circ \frac{\dd x_2}{x_2} \circ \dots \circ \frac{\dd x_n}{x_n}\right)\]
\[=(-1)^{n+1} \log \varphi\, \widehat{\re}_{n+1} \left( \int_\ii^1 \frac{\dd x_1}{x_1-\varphi}  \circ \frac{\dd x_2}{x_2} \circ \dots \circ \frac{\dd x_n}{x_n}  - \int_\ii^1\frac{\dd x_1}{x_1+ \varphi^{-1}}  \circ \frac{\dd x_2}{x_2} \circ \dots \circ \frac{\dd x_n}{x_n}\right) \] 
\begin{equation}\label{eq:golden1}
= (-1)^{n+1} \log \varphi \, \widehat{\re}_{n+1} ( \Li_n(\ii\varphi^{-1}) - \Li_n(\varphi^{-1})- \Li_n(-\ii\varphi) +\Li_n(-\varphi)). 
\end{equation}
We used the iterated integral notation in the previous equalities.

The last sum in (\ref{eq:golden})  yields
\[ \sum_{k=2,\, k \,\mathrm{even}}^{n}\frac{(-1)^{n-k+1}2^{k} B_{k} \log^k \varphi }{k!} \left( \int_\ii^1\re \left(\frac{\dd x_k}{x_k -\varphi}\right) \circ  \frac{\dd x_{k+1}}{x_{k+1}} \circ \dots \circ \frac{\dd x_{n}}{x_{n}}  \right.\]
\[ \left. + (-1)^k \int_\ii^1 \re \left(\frac{\dd x_k}{x_k+\varphi^{-1}} \right ) \circ \frac{\dd x_{k+1}}{x_{k+1}}\circ \dots \circ \frac{\dd x_{n}}{x_{n}}   \right)\]
\[=(-1)^{n+1}\sum_{k=2,\, k \,\mathrm{even}}^{n}\frac{2^{k} B_{k} \log^k \varphi }{k!} \widehat{\re}_{n-k+1}\left( \int_\ii^1 \frac{\dd x_k}{x_k -\varphi} \circ  \frac{\dd x_{k+1}}{x_{k+1}} \circ \dots \circ \frac{\dd x_{n}}{x_{n}}  \right.\]
\[ \left. + \int_\ii^1 \frac{\dd x_k}{x_k+\varphi^{-1}}  \circ \frac{\dd x_{k+1}}{x_{k+1}}\circ \dots \circ \frac{\dd x_{n}}{x_{n}}   \right)\]
\begin{equation}\label{eq:golden2}
 = (-1)^{n+1}\sum_{k=2,\, k \,\mathrm{even}}^{n}\frac{2^{k} B_{k} \log^k \varphi }{k!} \widehat{\re}_{n-k+1}(\Li_{n-k+1}(\ii\varphi^{-1}) - \Li_{n-k+1}(\varphi^{-1}) + \Li_{n-k+1}(-\ii\varphi) -\Li_{n-k+1}(-\varphi)       ).
\end{equation}

The sum of Eqs. (\ref{eq:golden1}) and (\ref{eq:golden2}) is 
\[= (-1)^{n+1}\sum_{k=1}^{n}\frac{2^{k} B_{k} \log^k \varphi }{k!} \widehat{\re}_{n+1}(\Li_{n-k+1}(-\ii \varphi)- \Li_{n-k+1}(-\varphi))\]
\[+(-1)^{n+1}\sum_{k=1}^{n}\frac{2^{k} B_{k}\log^k \varphi^{-1} }{k!} \widehat{\re}_{n+1}(\Li_{n-k+1}(\ii \varphi^{-1})- \Li_{n-k+1}(\varphi^{-1}))\]

Now the generalized Mahler measure is equal to  $n! \left( \frac{2}{\pi\ii}\right)^n$ times
\[-\frac{1}{2^n} \mathcal{L}_{n+1}(\varphi^2)+ \frac{1}{2^n}\sum_{k=3,\, k \,\mathrm{odd}}^{n+1} \frac{(\ii \pi)^{n-k+1}}{(n-k+1)!}\mathcal{L}_k(-\varphi^2)\]
\[+(-1)^{n+1}(\mathcal{L}_{n+1}(-\ii \varphi)- \mathcal{L}_{n+1}(-\varphi)+ \mathcal{L}_{n+1}(\ii \varphi^{-1})- \mathcal{L}_{n+1}(\varphi^{-1}))  \]
\[-(-1)^{n+1}\widehat{\re}_{n+1} (\Li_{n+1}(-\ii \varphi)- \Li_{n+1}(-\varphi)+ \Li_{n+1}(\ii \varphi^{-1})- \Li_{n+1}(\varphi^{-1})).\]

After  simplifying, we obtain:
\[ m\left(1+x_1-x_1^{-1}, \dots, 1+x_{n}-x_{n}^{-1}\right) = n! \sum_{k=3}^{n} \frac{1}{(n-k+1)! (\ii \pi)^{k-1}}\mathcal{L}_k(-\varphi^2)\]
\begin{equation}
+n! \left( \frac{2\ii}{\pi}\right)^n\widehat{\re}_{n+1} (\Li_{n+1}(-\ii \varphi)- \Li_{n+1}(-\varphi)+ \Li_{n+1}(\ii \varphi^{-1})- \Li_{n+1}(\varphi^{-1})). 
\end{equation}
A particular case is with $n=1$, where we obtain
\[ m(1+x_1-x_1^{-1})= -\log \varphi =\log (1+\varphi) ,\]
thus recovering the expected result from Mahler measure.

For $n=2$ we obtain
\[m\left(1+x_1-x_1^{-1},  1+x_{2}-x_{2}^{-1}\right) = \frac{2}{\pi^2}(\Li_3(\varphi^2) - \Li_3(-\varphi^2)) -\log \varphi,\]
which recovers the Mahler measure of $1+x_1-x_1^{-1}+z(1+x_{2}-x_{2}^{-1})$.
\subsection{A limit formula}

We end this section with an observation regarding the behavior of the generalized Mahler measure when the number of involved polynomials  goes to infinity.

\begin{prop} \label{prop} Let $P \in \Cset(x_1,\dots, x_n)$ and let $f_i=P(x_{i,1}, \dots, x_{i,n})$ for $i=1, \dots r$. Then
\[\lim_{r \rightarrow \infty} m(f_1,\dots,f_r) = \log ||P||_\infty\]
where $||\cdot||_\infty$ stands for the sup norm on $\TT^n$.
\end{prop}
\begin{pf}  For a measurable subset $S$  of $(\TT^n)^r$, let $I_S$ denote the integral over $S$:
\[ I_S = \frac{1}{(2 \pi \ii)^{nr}} \int_{S} \max_{1 \leq i \leq r} \log  |P(x_{i,1}, \dots, x_{i,n})|\frac{\dd x_{1,1}}{x_{1,1}} \dots \frac{\dd x_{r,n}}{x_{r,n}}.\]
Then we need to investigate $m(f_1,\dots,f_r)= I_{(\TT^n)^r}$ as $r$ goes to infinity.

First observe that
 \[m(f_1,\dots,f_r) \leq \log ||P||_\infty.\]
Then we just need to understand the lower bound. Consider the set 
\[ E_\epsilon := \{(x_1, \dots, x_n) \in \TT^n \, | \, \, |P(x_1, \dots, x_n)| \geq ||P||_\infty - \epsilon \}.\]  
Let $m_\epsilon$ denote the Lebesgue measure of $E_\epsilon$ (normalized so that the torus has measure one).

We write  $(\TT^n)^r$ as a union of sets: 
\begin{eqnarray*}
(\TT^n)^r &=& E_\epsilon \times (\TT^n)^{r-1} \cup (\TT^n \setminus E_\epsilon) \times E_\epsilon \times (\TT^n)^{r-2} \cup \dots \cup (\TT^n\setminus E_\epsilon)^r, \\
&=&A \cup (\TT^n\setminus E_\epsilon)^r.
\end{eqnarray*}

Notice that the set $A$ has measure 
\[m_\epsilon + (1- m_\epsilon)m_\epsilon+ (1-m_\epsilon)^2 m_\epsilon + \dots + (1-m_\epsilon)^{r-1} m_\epsilon = 1-(1-m_\epsilon)^r.\]

Therefore,
\begin{equation}\label{ineq3}
I_A \geq (1-(1-m_\epsilon)^r) \log ( ||P||_\infty - \epsilon).
\end{equation}
 
Outside the set $E_\epsilon$ we do not have control of the value of $\log|P|$. This is a problem  if the polynomial has zeros in the torus. Let
\[Z_\eta = \{ (x_1, \dots, x_n) \in \TT^n \, | \, \, |P(x_1, \dots, x_n)| \leq \eta \},\]
and let $m_\eta$  denote the measure of $Z_\eta$.

We write
\[ (\TT^n\setminus E_\epsilon)^r = (\TT^n\setminus E_\epsilon)^r \setminus Z^r_\eta \cup Z^r_\eta.\]
Now 
\begin{equation} \label{ineq2}
I_{(\TT^n\setminus E_\epsilon)^r \setminus Z^r_\eta}\geq ((1-m_\epsilon)^r-m_\eta^r) \log \eta .
\end{equation}

On the other hand,
\begin{eqnarray}\label{ineq}
 I_{Z^r_\eta}  &\geq&  \frac{1}{(2 \pi \ii)^{rn}} \int_{Z_\eta^r} \log  |P(x_{1,1}, \dots, x_{1,n})|\frac{\dd x_{1,1}}{x_{1,1}} \dots \frac{\dd x_{r,n}}{x_{r,n}} \nonumber\\
 & =& m_\eta^{r-1} \frac{1}{(2 \pi \ii)^{n}} \int_{Z_\eta} \log  |P(x_{1}, \dots, x_{n})|\frac{\dd x_{1}}{x_{1}} \dots \frac{\dd x_{n}}{x_{n}} . 
\end{eqnarray}

First choose a small $\epsilon$ and large $r$, so $I_A$ becomes close to $\log ||P||_\infty$ due to inequality (\ref{ineq3}).

Now fix a small $\eta$. We have that both $(1-m_\epsilon)^r$ and $m_\eta^r$ are small for $r$ large. If $r$ is large enough,  the right side of inequality (\ref{ineq2}) is close to zero. Then either $I_{(\TT^n\setminus E_\epsilon)^r \setminus Z^r_\eta}$ is positive or has a small absolute value.

Finally, $m_\eta^{r-1}$ is small for $r$ large and the right side of  inequality (\ref{ineq}) is close to zero  because $\log |P| \in L^1(\TT^n)$.  Then either $ I_{Z^r_\eta}$ is positive or has a small absolute value. 

We reach our conclusion by considering $ I_{(\TT^n)^r} = I_A+I_{(\TT^n\setminus E_\epsilon)^r \setminus Z^r_\eta} +  I_{Z^r_\eta} $.

\end{pf}

We may conclude interesting results. For instance, from Eq. (\ref{limit}) we deduce
\begin{equation}
\lim_{m \rightarrow \infty} \sum_{j=1}^{m-1} (-1)^{j}\binom{2m-1}{2j} \frac{(2j)!(1-2^{2j})}{(2\pi)^{2j}} \zeta(2j+1) = \log 2.
\end{equation}

\bigskip
\begin{ack}

I thank my Ph.D supervisor Fernando Rodriguez-Villegas for many incredibly helpful and encouraging discussions. I would also like to acknowledge enlightening discussions with David Boyd, Herbert Gangl and Vincent Maillot about Mahler measure, Bloch groups, and cohomology, and with Enrico Bombieri about generalized Mahler measure and more specifically, Proposition \ref{prop}. I am very grateful to the referee whose dedicated and careful review has much contributed to the clarity of this work (including necessary corrections), and for the proof of Eq. (\ref{eq:22}).

Most of this work is part of my Ph. D. dissertation at the Department of Mathematics at
the University of Texas at Austin. I am indebted to John Tate and the Harrington
fellowship for their generous financial support during my graduate studies.

Part of this research was completed while I was a member at the Institute for Advanced Study, a visitor at the Institut des Hautes \'Etudes Scientifiques, and a Postdoctoral Fellow at the Pacific Institute for the Mathematical Sciences. I am grateful for their support. 

\end{ack}


\begin{thebibliography}{XXXXXX}



\bibitem[Boy81]{B1} D. W. Boyd, Speculations concerning the range of Mahler's measure, {\em Canad.  Math. Bull.\/} {\bf 24} (1981), 453--469.

\bibitem[Boy98]{B2} D. W. Boyd, Mahler's measure and special values of  L-functions, {\em Experiment. Math.\/} {\bf 7} (1998), 37--82.

\bibitem[BR-V02]{BRV1} D. W. Boyd, F. Rodriguez-Villegas, Mahler's measure and the dilogarithm (I), {\em Canad. J. Math.\/} {\bf 54} (2002), 468--492.

\bibitem[BR-V03]{BRV2} D. W. Boyd, F. Rodriguez-Villegas, with an appendix by N. M. Dunfield, Mahler's measure and the dilogarithm (II), (preprint, July 2003).


\bibitem[Con03]{C} J. Condon, Calculation of the Mahler measure of a three variable polynomial, (preprint, October 2003).



\bibitem[Den97]{D} C. Deninger,  Deligne periods of mixed motives, $K$-theory and the entropy of certain $ Z\sp n$-actions, {\em  J. Amer. Math. Soc. \/} {\bf 10}  (1997),  no. 2, 259--281.


\bibitem[GO04]{GO} Y. Gon, H. Oyanagi, Generalized Mahler measures and multiple
sine functions. {\em Internat. J. Math.\/} {\bf 15} (2004), no. 5, 425--442.



\bibitem[Gon95]{G3} A. B. Goncharov, Geometry of Configurations, Polylogarithms, and Motivic Cohomology, {\em Adv. Math.\/}  {\bf 114}  (1995),  no. 2, 197--318.



\bibitem[Gon02]{G4} A. B. Goncharov, Explicit regulator maps on polylogarithmic motivic complexes.  {\em Motives, polylogarithms and Hodge theory, Part I (Irvine, CA, 1998),\/}  245--276, {\em Int. Press Lect. Ser., 3, I, Int. Press, Somerville, MA,} 2002.


\bibitem[Gon05]{G5} A. B. Goncharov, Regulators. {\em Handbook of $K$-theory}. Vol. 1, 2,  295--349, Springer, Berlin, 2005.


\bibitem[Lal03]{L} M. N. Lal\'{\i}n, Some examples of Mahler measures as multiple polylogarithms, {\em  J. Number Theory\/} {\bf 103}  (2003),  no. 1, 85--108.


\bibitem[Lal06]{L4} M. N. Lal\'{\i}n, On certain combination of colored multizeta values, {\em J. Ramanujan Math. Soc.\/} {\bf 20} (2006), no. 1, 115--127.

\bibitem[Lal07]{L3} M. N. Lal\'{\i}n, An algebraic integration for Mahler 
measure, (to appear in {\em Duke Math. J.} {\bf 138} (2007))



\bibitem[R-V97]{RV} F. Rodriguez-Villegas, Modular Mahler measures I, {\em Topics in number theory (University Park, PA 1997),\/} 17--48, {\em Math. Appl., 467, Kluwer Acad. Publ. Dordrecht,\/} 1999.



\bibitem[Smy81]{S1} C. J. Smyth, On measures of polynomials in several  variables, {\em Bull. Austral. Math. Soc. Ser. A\/} {\bf 23} (1981), 49--63. Corrigendum (with G. Myerson): {\em  Bull. Austral. Math. Soc.\/} {\bf 26} (1982), 317--319.
  
\bibitem[Smy02]{S2} C. J. Smyth, An explicit formula for the Mahler measure  of a family  of 3-variable polynomials, {\em J. Th. Nombres Bordeaux\/} {\bf 14} (2002), 683--700.

 








\bibitem[Zag91]{Z2} D. Zagier, Polylogarithms, Dedekind Zeta functions, and the Algebraic $K$-theory of Fields,  {\em Arithmetic algebraic geometry\/} (Texel, 1989),  391--430, Progr. Math., 89, Birkh\"auser Boston, Boston, MA, 1991.
 

\bibitem[Zha03]{Zh} J. Zhao, Goncharov's relations in Bloch's higher Chow group $CH^3(F,5)$, {\em  J. Number Theory\/} {\bf 124}  (2007),  no. 1, 1--25.

\end{thebibliography}
\end{document}